\documentclass[12pt]{article}
\usepackage[english]{babel}
\usepackage{amsmath,amssymb,latexsym}
\usepackage{a4format,arthead}
\usepackage{theorem}
\usepackage{humath}
\usepackage{local}
\usepackage{eulercal}
\begin{document}
\thispagestyle{empty}
\title{Traces on algebras of parameter dependent pseudodifferential
operators and the eta--invariant}

\author{Matthias Lesch \ and     Markus J.~Pflaum}
\date{10-21-1998\\{\small revised version}}
%
% version: 1.0, 29-04-1998
%

%%% Local Variables: 
%%% mode: latex
%%% TeX-master: "LesPflTAPDPOEI"
%%% End: 

%
\maketitle
\begin{abstract}
We identify Melrose's suspended algebra of pseudodifferential
operators with a subalgebra of
the algebra of parametric pseudodifferential operators with
parameter space $\R$.
For a general algebra of parametric pseudodifferential
operators, where the parameter space may now be a cone  $\Gamma\subset\R^p$,
we construct a unique ``symbol valued trace'',
which extends the $L^2$--trace on operators of small order.
This construction is in the spirit of a trace due to 
Kontsevich and Vishik in the nonparametric case. 
Our trace allows to construct various trace functionals in a systematic
way. Furthermore we study the higher--dimensional
eta--invariants on algebras with parameter space $\R^{2k-1}$.
Using Clifford representations we construct for each first order
elliptic differential operator a natural
family of parametric pseudodifferential operators over
$\R^{2k-1}$. The eta--invariant of this family
coincides with the spectral eta--invariant of the operator. 

%{\it Keywords and phrases.} 
%Pseudodifferential operators, noncommutative residue

\end{abstract}

%%% Local Variables: 
%%% mode: latex
%%% TeX-master: "LesPflTAPDPOEI"
%%% End: 

%
\tableofcontents
\newpage
\section{Introduction}

Let $M$ be a smooth compact Riemannian manifold without
boundary. Furthermore, let $E$ be a hermitian vector bundle
over $M$. We denote by $\CL^*(M,E)$ the algebra of classical
pseudodifferential operators acting on $L^2(M,E)$.
It is well--known that up to a scalar factor $\CL^*(M,E)$
has a unique trace, the residue trace of {\sc Guillemin} \cite{Gui:RTCAFIO}
and {\sc Wodzicki} \cite{Wod:NCRIFK}. The residue trace vanishes on trace class
operators and hence there is no trace on $\CL^*(M,E)$ which extends the
$L^2$-trace. However, {\sc Kontsevich} and {\sc Vishik} 
\cite{KonVis:DEPO,KonVis:GDEO} constructed a functional
$$
 \TR : \bigcup_{\alpha \in \R \setminus \Z} \CL^\alpha (M,E) \rightarrow \C ,
$$
which extends the $L^2$-trace and satisfies $\TR (AB) = \TR (BA)$ for 
$A,B \in \CL^* (M,E)$ with ${\rm ord} (A) + {\rm ord} (B) \notin \Z$.
$\TR$ has further interesting properties 
(cf.~\cite{KonVis:DEPO,KonVis:GDEO}, \cite[Sec.5]{Les:NRPOPS}). 

In this paper we study traces on the algebra of parameter
dependent pseudodifferential operators $\CL^*(M,E,\Gamma)$, where
$\Gamma\subset\R^p$ is a conic set.
These algebras play an important role in the study of the resolvent
of an elliptic differential operator in which case $\Gamma$
is a sector in $\C$ (cf. \cite{Shu:POST}). 

Our first result shows that $\CL^*(M,E,\R)$ contains a canonical isomorphic
image of the algebra $\CL^*_{\rm sus}(M,E)$ introduced by
{\sc Melrose} \cite{Mel:EIFPO}.
This algebra appears naturally in an index theorem for
manifolds with boundary \cite[Sec.~12]{MelNis:HPOIMB}.
It should be thought of a pseudodifferential suspension of the 
algebra $\CL^*(M,E)$.

More generally we then study the algebra $\CL^*(M,E,\Gamma)$ for
a connected cone $\Gamma\subset \R^p$ with nonempty interior.
Our main result is that for $H_{\rm dR}^1(\Gamma)\neq 0$
the algebra $\CL^*(M,E,\Gamma)$ 
has a unique ``symbol valued trace''. More precisely,
there exists a unique linear map 
$$\TR:{\rm CL}^*(M,E,\gG)\rightarrow {\rm PS}^*(\gG)/\C[\mu_1,...,\mu_p]$$
with the following properties:
\begin{enumerate}\renewcommand{\labelenumi}{{\rm (\roman{enumi})}}
\item $\TR(AB)=\TR(BA)$, i.e. $\TR$ is a ``trace'',
\item $\TR(\pl_j A)=\pl_j \TR(A)$, \quad $j=1,...,p$,
\item If $A\in {\rm CL}^m(M,E,\gG)$ and $m+\dim M<0$ then
$$\TR(A)(\mu)=\trltwo(A(\mu)).$$
\end{enumerate}
Here ${\rm PS}^*(\gG)$ is the class of symbols having a complete
asymptotic expansion in terms of homogeneous functions and $\log$--powers.

So the parametric situation is different from the nonparametric one:
the symbol valued $L^2$--trace can be extended, however only
modulo polynomials.
$\TR$ should be viewed as the analogue of the Kontsevich-Vishik trace, since
it can be constructed quite similarly (see Remark \ref{ML5-4.4} below).
Our main result allows to construct various traces on the
algebra $\CL^*(M,E,\Gamma)$ just by composing $\TR$ with linear functionals
on  ${\rm PS}^*(\gG)/\C[\mu_1,...,\mu_p]$.

The most important examples are the extended and the formal trace
$\overline{\Tr}$ resp. $\widetilde{\Tr}$. For $A\in \CL^m(M,E,\R^p)$
the extended trace is given by
\begin{equation}
   \overline{\Tr}(A)=\regint_{\R^p} \TR(A)(\mu) d\mu,
\end{equation}
where $\reginttext_{\R^p}$ is a certain regularization of the 
integral (cf.~\ref{S2-4.3}).
If $m+\dim M+p<0$ the function $\TR(A)$ is integrable, and 

\begin{equation}
     \overline{\Tr}(A)=\int_{\R^p} \trltwo(A(\mu)) d\mu.
\end{equation}
holds indeed.

From $\CL^*(M,E,\R^p)$ we can construct a de Rham complex
$(\Omega^*\CL^*(M,E,\R^p),d)$ in a canonical way. Then
$\overline{\Tr}$ extends to a graded trace on the complex 
$\Omega^*\CL^*(M,E,\R^p)$
\begin{equation}
  \overline{\Tr} (\go):=\casetwo{0}{\mbox{if } \: {\rm deg }\,\go <p,}{
    \overline{\Tr}(f)}{\mbox{if }\: \go=fdx_1\wedge\ldots\wedge dx_p.}
\end{equation}
However, the graded trace is not closed, but its derivative
$\widetilde\Tr:=d\overline{\Tr}$ is a closed graded
trace on $\Omega^*\CL^*(M,E,\R^p)$.
If $p=1$, then $\overline{\Tr}$ and $\frac{1}{2\pi}\widetilde\Tr$
coincide with the corresponding traces introduced in
\cite{Mel:EIFPO}. 

Like in \cite{Mel:EIFPO} $\widetilde\Tr$ is an analogue of the residue trace.
It only depends on finitely many terms of the symbol expansion of the
operator. One of the results of {\sc Melrose} \cite{Mel:EIFPO}
was the construction of the eta--homomorphism. In our notation
\begin{equation}
     \eta:\CL^*(M,E,\R)^{-1}\longrightarrow \C
\end{equation}
is a homomorphism from the group of invertible elements of
$\CL^*(M,E,\R)$ into the additive group $\C$.
In some sense $\eta$ generalizes the winding number. Namely,
for $A\in \CL^*(M,E,\R)^{-1}$ one has
\begin{equation}
\begin{split}
      \eta(A)=&\frac{1}{\pi i}  \overline\Tr(A^{-1}dA)
             = \frac{1}{\pi i}  \regint_\R \TR(A^{-1} A')(x) dx.
\end{split}
\end{equation}
In case $A$ is a function on $\R$ taking values in the space of invertible 
matrices and which is constant outside a compact set then $\frac 12 \eta(A)$ 
is an integer equal to the winding number of $A$.
Thus it is natural to expect a similar invariant for odd--dimensional
parameter spaces. Indeed for $A\in \CL^*(M,E,\R^{2k-1})^{-1}$
we put
\begin{equation}
    \eta_k(A):= 2 c_k \overline{\Tr}((A^{-1}dA)^{2k-1}),
    \label{ML5-1.5}
\end{equation}
where $c_k$ is a normalization constant. Again, if $A$ is
just a matrix valued function and constant outside a compact set,
\myref{ML5-1.5} is an even integer which actually classifies
the $(2k-1)$th homotopy group of ${\rm GL}(\infty,\C)$.

In contrast to its finite--dimensional analogue $\eta_k$ is not
a homotopy invariant. However its variation is local that means for 
a smooth family $A_s$ of invertible elements the equality
 \begin{equation}
    \frac{d}{ds} \eta_k (A_s) = 2(2k-1) c_k \widetilde{\Tr} \left(
    (A_s^{-1} \partial_s A_s) (A_s^{-1} d A_s)^{2k-2} \right)
  \end{equation}
holds true.
Unfortunately $\eta_k$ is not a homomorphism for $k\ge 2$, instead we have
\begin{equation}
     \eta_k(AB)-\eta_k(A)-\eta_k(B)=\widetilde\Tr(\go(A,B)),
\end{equation}
where $\go(A,B)$ denotes a universal polynomial in the $1$--forms
$B^{-1} (A^{-1}dA)B, B^{-1}dB$.
So the defect of the additivity is a symbolic term.

Finally we compare $\eta_k$ with the spectral eta--invariant.
For any first order invertible self--adjoint elliptic differential
operator $D$ we construct a natural family
$\cD(\mu):=D+c(\mu)$ in $\CL^1(M,E,\R^{2k-1})$ such that
\begin{equation}
     \eta_k(\cD)=-\eta(D),
\end{equation}
where $c$ is the standard Clifford representation and $\eta(D)$
the spectral eta--invariant of $D$.

\bigskip
We understand that some of our results also have been obtained by
{\sc R. B. Melrose} and {\sc V. Nistor} \cite{MelNis:IP}.

\bigskip \noindent
{\bf Acknowledgement:} 
The first named author gratefully acknowledges the hospitality
and financial support of the Erwin--Schr\"odinger Institute, Vienna, 
where part of this work was completed. Furthermore, the first named author
was supported by Deutsche Forschungsgemeinschaft. 

%%% Local Variables: 
%%% mode: latex
%%% TeX-master: "LesPflTAPDPOEI"
%%% End: 

%
\section{Review of parametric pseudodifferential operators}
The concept of parameter dependent symbols and pseudodifferential
operators used in this article involves several different classes
of symbol spaces. For the convenience of the reader and to fix the
notation we briefly recall some basic facts about symbols
and the corresponding operator calculus. As  general references we mention
the books {\sc Shubin} \cite{Shu:POST} and {\sc Grigis-Sj{\o}strand}
\cite{GriSjo:MADO}. \vspace{2mm}

A conic manifold is a smooth principal fiber bundle $\Gamma \rightarrow B$
with structure group $\R_+:=(0,\infty)$. 
It is always trivializable (cf.~{\sc Duistermaat}
\cite{Dui:FIO}, \S 2.1). A subset $\Gamma \subset \dot{\R}^{\nu} := \R^{\nu}
\setminus \{ 0 \} $ which is a conic manifold by the natural $\R_+$-action 
on $\dot{\R}^{\nu}$ is called a conic set. 
The base manifold of a conic set $\Gamma \subset \dot{\R}^{\nu}$ is
isomorphic to $S \Gamma := \Gamma \cap S^{\nu -1}$. By a cone 
$\Gamma \subset \R^{\nu}$ we will always mean a conic set or the 
closure of a conic set in $\R^{\nu}$ such that $\Gamma$ has nonempty interior.
Thus $\R^n$ and $\dot{\R}^n$ are cones, but only the latter is a conic set. 

Now let $M$ be a smooth manifold, $m \in \R$, $\Gamma \subset \R^{\nu}$ a 
cone, and $0 < \rho \leq 1$. Then by 
$\sym^m_{\gr} (X,\gG)$ we denote the space of all functions
$a (x,\xi) \in C^{\infty} (M \times \gG)$ such that 
for every differential operator $D$ on $M$, all compact $L\subset \Gamma$
and $K \subset M$ we have the uniform estimate
\begin{equation}
  \left| D \, \partial^{\ga}_{\xi} a (x, \xi) \right| \: \leq \:
  C_{K,L,D,\ga} \, \brac{\xi}^{m-\gr |\ga|} , \quad x \in K, \: \xi \in 
  L^{\rm c}, \: \alpha \in \Z_+^{\nu}.
\end{equation}
Here, $L^{\rm c} := 
\{ t \xi \, | \, \xi \in L, t \geq 1 \}$, 
$C_{K,L,D,\alpha} >0$ and $\brac{\xi} := ( 1 + |\xi |^2 )^{1 / 2}$ for all 
$\xi \in \R^{\nu}$. In case $\gr =1$ we write $\sym^m (M, \Gamma)$ 
for $\sym^m_{\gr} (M,\Gamma)$.

A symbol $a \in \sym^m (M, \Gamma)$ is called classical polyhomogeneous of 
degree $m$ or just classical, if it has an asymptotic expansion of the form
$  a \sim \sum_{j\ge 0} \, a_{m-j} $, where the $a_k (x,\xi) \in \sym^k
(M,\Gamma)$ are $k$-homogeneous in $\xi$ of degree $k$.
The space of classical polyhomogeneous symbols of order $m$ is denoted by
$\csym^m (M,\Gamma)$. 

Now let $U \subset \R^n$ be an open set, and 
$a \in \sym^m_{\rho} (U , \R^n \times \Gamma)$.
For each fixed $\mu_0$ we have $a(\cdot, \cdot, \mu_0) \in \sym^m (U, \R^n)$, 
hence we obtain a family of pseudodifferential operators
parametrized over $\Gamma$ by putting
\begin{equation}
\begin{split}
  \big[ \Op( a(\mu) ) \, u \big] \, (x) & :=  \big[ A(\mu) \, u \big] (x)
  := \int_{\R^n} \, \e^{i \langle x,\xi \rangle} \, a(x,\xi,\mu) \, \hat{u} (\xi ) 
  \, \dbar \xi , \quad \dbar \xi  := (2 \pi)^{-n} d \xi .       
\end{split}
\end{equation}
Note that for $a \in \sym^{-\infty} (U , \R^n \times \Gamma)$ the
operator $A(\mu)$ has
a kernel in $\cS (\Gamma , C^{\infty} (U \times U))$, the Schwartz space
of $C^{\infty} (U \times U)$-valued functions. We denote by 
${\rm L}^m_{\rho} (U, \Gamma)$ the set of all $A(\mu)$, where 
$a \in \sym^m_{\rho} (U, \R^n \times \Gamma)$. 
In case $\Gamma = \{ 0 \}$ we obtain the well-known space 
${\rm L}^m_{\rho} (U)$  of pseudodifferential operators of order $m$ 
and type $\rho$ on $U \subset \R^n$.

For a smooth manifold $M$ and vector bundles $E, F$ over $M$ the spaces
${\rm L}^m (M,E,F; \Gamma)$ of parameter dependent pseudodifferential
operators between sections of $E,F$ are defined in the usual way
by patching together local data.

The space of parameter dependent pseudodifferential operators with symbol
lying in the space $\csym^m (U,\R^n \times \Gamma)$ will be denoted by
${\rm CL}^m (U, \Gamma)$. Its elements are the classical parameter dependent
pseudodifferential operators over $U$. 
Following {\sc Grubb} and {\sc Seeley} \cite{GruSee:WPPOAPSBP} we also call
these operators strongly polyhomogeneous.
% (later we will also consider weakly polyhomogeneous operators).
\begin{example}
  Let $M$ be a compact manifold,
  $A \in {\rm CL}^m (M)$, and assume that $\Gamma
  \subset \C \setminus \{ 0 \} $ is a cone such that $\sigma^m_A (x,\xi) -z $
  is invertible for $z \in \Gamma$. If $A$ is a differential operator and
  ${\rm spec} \, A \cap \Gamma = \emptyset$ then $(A -z )^{-1} \in {\rm L}^{-m}
  (M , \Gamma)$. However, in general this need not be true for 
  pseudodifferential operators.
\end{example}
The following result is just a mild generalization of the classical
resolvent expansion of a differential operator
(see e.g. \cite[Sec.~1.7]{Gil:ITHEASIT}).
\begin{theorem} \label{S1-1.2}
  Let $M$ be a compact manifold,
${\rm dim}\, M=:n$, and $A \in {\rm L}^m (M,E,\Gamma)$ (resp.
${\rm CL}^m (M,E,\Gamma)$). If $m +
  {\rm dim} \, M < 0$ then $A(\mu)$ is trace class for all $\mu \in \Gamma$
  and
  \begin{displaymath}
     \tr \, A (\cdot) \in {\rm S}^{m + {\rm dim} \, M} (\Gamma)\quad
(\mbox{\rm resp.}\: {\rm CS}^{m + {\rm dim} \, M} (\Gamma)).
  \end{displaymath}
\end{theorem}
\begin{proof} 
 We present 
 the proof for $A\in {\rm CL}$. For $A\in {\rm L}$ it is
 even a bit simpler.
 Choosing a suitable partition of unity it suffices to prove the claim
 for $E = \C$, $M =U$ a coordinate patch, and $A$ compactly supported,
 i.e.
 \begin{equation}
    (A u) (x) = \int_{\R^n} \, \sigma_A (x,\xi,\mu) \, \hat{u}(\xi)
    \, \dbar \xi
    = \int_{\R^n} \int_{U} \, {\rm e}^{i \langle x-y,\xi \rangle} \, \sigma_A 
    (x, \xi, \mu) \, u(y) \, dy \, \dbar \xi,
  \label{ML5-1.3}  
 \end{equation}
  where $\sigma_A \in {\rm CS}^m (U, \R^n \times \Gamma)$ and
  $\pi_1 ( {\rm supp} \, \sigma_A ( \cdot, - , \mu) )  \subset K \subset U$
  is compact for every $\mu \in \Gamma$. 
  For fixed $\mu$ we have $A(\mu) \in {\rm CL}^m (U)$, hence $A(\mu)$ is
  trace class since $m < -{\rm dim} \, M = -n$.
  Since $\sigma_A \in {\rm CS}^m (U, \R^n \times \Gamma)$ we have
  \begin{equation}
    \sigma_A \sim \sum_{j=0}^\infty \, a_{m-j}
  \end{equation}
  with $a_{m-j} (x ,\lambda \xi , \lambda \mu ) = \lambda^{m-j} \, a_{m-j}
  (x , \xi , \mu )$ for $\lambda \geq 1$, $ | ( \xi, \mu )| \geq 1$.
  Thus we write
  \begin{equation}
    \sigma_A = \sum_{j=0}^{N-1} \, a_{m-j} \, + R_N
  \end{equation}
  with $R_N \in {\rm CS}^{m-N} (U , \R^n \times \Gamma )$.
  Now pick $L \subset \Gamma$ compact and a multiindex $\alpha$. Then 
  \begin{equation}
  \begin{split}
    \left| \partial^{\alpha}_{\mu} \tr \, {\rm Op} (R_N (\mu)) \right|
    & = \left| \int_K \int_{\R^n} \, \partial_{\mu}^{\alpha} R_N
    (x , \xi , \mu ) \dbar \xi \, dx \right| \\
    & \leq C_{\alpha, K , L} \int_K \int_{\R^n} ( 1 + ( |\xi|^2
      +|\mu|^2)^{1/2} )^{m -|\alpha| - N} \, \dbar\xi \, dx \\
    & \leq C_{\alpha,K,L} (1 + |\mu|)^{m+n - |\alpha| - N }.
  \end{split}\label{G1-1.7}
  \end{equation}
  Furthermore let $\lambda \geq 1$, $| \mu  | \geq 1$. Then
  \begin{equation}
  \begin{split}
    {\rm tr} \, {\rm Op} (a_{m-j} (\lambda \mu)) & = \int_U \int_{\R^n}
    \, a_{m-j} (x, \xi, \lambda \mu ) \, \dbar \xi \, dx   \\
    & = \lambda^{m-j} \, \int_U \int_{\R^n} \, a_{m-j} (x , \lambda^{-1} \xi ,
    \mu ) \, \dbar \xi \, dx \\
    & = \lambda^{m+n-j} \, {\rm tr} \, {\rm Op} (a_{m-j} (\mu )) ,
  \end{split}
  \end{equation}
and similar to \myref{G1-1.7} one shows that $\tr({\rm Op}(a_{m-j}))\in
\sym^{m-j+n}(\gG)$, 
  thus
  \begin{equation}
    {\rm tr} \, A(\mu) \sim \sum_{j=0}^\infty \, {\rm tr} \, {\rm Op}
    (a_{m-j} (\mu) ) 
  \end{equation}
  in ${\rm CS}^{m+n} (\Gamma)$, where ${\rm tr} \, {\rm Op} (a_{m-j} (\mu) )$
  is homogeneous of degree $m+n-j$ for $\mu \geq 1$.
\end{proof}
The previous proof provides even more, namely
\begin{theorem} \label{S1-1.3}
  Let $M$ be a smooth manifold. If $m + {\rm dim} \, M < 0$ then for 
  any properly supported $A \in {\rm L}^m (M,E,F;\Gamma)$ (resp. 
${\rm CL}^m (M,E,F;\Gamma)$) there is a 
  density 
  \begin{displaymath}\begin{split}
    &\omega_A \in {\rm S}^m (M, {\rm Hom} (E,F) \otimes | \Omega M | ; \Gamma )
       \\
    &(\mbox{\rm resp.}\:
  {\rm CS}^m (M, {\rm Hom} (E,F) \otimes | \Omega M | ; \Gamma ))
                  \end{split}
\end{displaymath}
  with the following properties:
  \begin{enumerate}\renewcommand{\labelenumi}{{\rm \arabic{enumi}.}}
    \item
      For any local chart $\phi : M \supset U \rightarrow \phi (U) \subset
      \R^n$ we have
    %  \begin{equation}
     $$   \phi^* \omega_{\phi_* A} = \go_A.
     $$% \end{equation}
    \item
      If $E = F$, then 
     % \begin{equation}
      $$  {\rm tr} \, A( \mu) = \int_M \, {\rm tr}_{E_x} \, 
        \omega_A (x, \mu) .
     $$% \end{equation}
  \end{enumerate}
\end{theorem}
\begin{proof}
  The proof of Theorem \ref{S1-1.2} shows that we can put
%  \begin{equation}
 $$   \omega_A (x, \mu) = \int_{\R^n} \, 
    \sigma_A (x , \xi,\mu ) \, \dbar\xi\,|dx|.
$$%  \end{equation}
  Then 1. and 2. follow easily.
\end{proof}
\begin{remark}
  For the preceding two theorems the assumption $m + {\rm dim} \, M < 0$
  was essential. However in this paper we will show that these theorems
  can be extended to arbitrary 
%(weakly) 
  parametric operators.
\end{remark}

%%% Local Variables: 
%%% mode: latex
%%% TeX-master: "LesPflTAPDPOEI"
%%% End: 

%
\section{Melrose's suspended algebra of pseudodifferential operators}
In the paper \cite{Mel:EIFPO} {\sc R. B. Melrose}
invented a "suspended" algebra
of \psdo s on a compact manifold. He introduced trace functionals on this
algebra and constructed the "$\eta$--homomorphism". In this section we 
will briefly recall the definition of the suspended algebra and
we will show that it is isomorphic to a subalgebra of $\CL^*(M,\R)$.

In the subsequent sections we will construct trace functionals on 
$\CL^*(M,\Gamma)$ which generalize the {\sc Melrose} traces.

Let $M$ be a compact manifold. Following \cite{Mel:EIFPO} 
$\Lsus^m(M)$ consists of those operators $A\in \CL^m(M\times\R)$ such
that
\renewcommand{\labelenumi}{{\rm (\roman{enumi})}}
\begin{enumerate}
\item $A$ acts as convolution in the second variable, i.e. by slight
abuse of notation
  \begin{equation}
     (Au)(x,t)=\int_M \int_\R K_A(x,y,t-s) u(y,s) \, ds dy,
    \label{ML1-G1.1}
\end{equation}
where $K_A$ denotes the convolution kernel of A.
\item The kernel satisfies
\begin{equation}   
   K_A\in C_0^\infty (M \times M \times \R; |\Omega M| \boxtimes 1)'
   +\cS(M^2\times \R; 1 \boxtimes |\Omega M | ),
   \label{ML1-G1.2}
\end{equation}
where $C_0^\infty (M \times M \times \R; |\Omega M| \boxtimes 1)' $ 
denotes the space of distributions which act on the smooth compactly supported
sections of the exterior tensor product of the density bundle $|\Omega M|$
with the trivial line bundle over $M \times \R$.
\end{enumerate}
For vector bundles $E$ and $F$, $\Lsus^m(M,E,F)$ is defined accordingly.
$\Lsus^*(M,E)$ is an order filtered algebra.
For $A\in \Lsus^*(M,E)$ {\sc Melrose} introduced what
he calls the indicial family $\widehat A$. This is the partial
Fourier transform in the $t$--variable. Namely, for $\mu\in\R$
we obtain a \psdo\ $\widehat A(\mu)\in \CL^m(M,E)$ by putting
\begin{equation}
\widehat A(\mu) g= e^{-it\mu} A(e^{i\mu(\cdot)} g)(\cdot,t),\quad
   g\in C^\infty(M,E).
\end{equation}

\begin{prop}\label{ML1-S1.1} The map
$$\Lsus^*(M,E)\longrightarrow \CL^*(M,E,\R),\quad A\mapsto \widehat A(\cdot)$$
is an order preserving injective homomorphism of $*$-algebras.
\end{prop}
\begin{proof} We only have to check that for $A\in \Lsus^m(M,E)$ we have
$\widehat A(\cdot)\in \CL^m(M,E,\R)$. It suffices
to check this locally for $E=\C$. Let $U$ be a coordinate patch in $M$.
Then for $u\in C^\infty_0(U\times\R)$ we have
\begin{equation}
     (Au)(x,t)=\int_{U\times\R} \sigma_A((x,t),(\xi,\mu))\hat u(\xi,\mu)
    e^{i(\langle x,\xi \rangle +t\mu)} \dbar \mu \dbar \xi,
  \label{ML1-G1.3}
\end{equation}
with $\sigma_A\in \CS^m(U\times\R, \R^n\times\R)$. In view of 
\myref{ML1-G1.1}, $\sigma_A$ is independent of $t$, hence
$\sigma_A(x,\xi,\mu)\in \CS^m(U,\R^n\times\R)$ and this is the
complete symbol of $\widehat A(\mu)$.
\end{proof}

Summing up the suspended algebra can be viewed as an algebra
of strongly polyhomogeneous parameter dependent \psdo s. 
We did not try to express (ii) in terms of the indicial family
$\widehat A$.
However, it turns out that the extended trace and the
$\eta$--homomorphism can be constructed without using
(ii). So one could equally well consider the algebra $\CL^*(M,E,\R)$
as $\Lsus^*(M,E)$.

There is no reason to restrict the consideration to $\R$ as a parameter
space. Therefore, we will deal with
$\CL^*(M,\Gamma)$ in the sequel .

%%% Local Variables: 
%%% mode: latex
%%% TeX-master: "LesPflTAPDPOEI"
%%% End: 

%
\section{Tracial maps on $\bf {\rm \bf L}^*(M,E,\Gamma)$}
During the whole section let $\Gamma\subset \R^p$ be a {\em connected}
cone with nonempty interior.

\begin{dfn} \label{S1-2.0} 
We define the following spaces of functions on $\Gamma$:
\begin{eqnarray*}
&&\cP^m=\cP^m(\gG):=\{P\restr \gG\,|\,
P\in\C[x_1,\ldots,x_p], {\rm deg}\, P\le m\},\quad m\in\R_+,\\
&&\tilde\sym^m(\Gamma):=\casetwo{\sym^m(\Gamma)}{m\not\in\Z_+,}{
  \{f\in \mcap\limits_{\eps>0} \sym^{m+\eps}(\Gamma)\,|\, \partial^\ga f\in
  \sym^{m-|\ga|}(\Gamma),|\ga|\ge m+1\}}{m\in\Z_+.}
\end{eqnarray*}
\end{dfn}
Note that $\tilde\sym^*(\Gamma)=\mcup\limits_{m\in\R}\tilde\sym^m(\Gamma)=
\sym^*(\Gamma)$.

Since $\Gamma$ is assumed to be connected with nonempty interior
the restriction map $C^\infty(\R^p)\to C^\infty(\Gamma), f\mapsto f\restr\gG$
induces an isomorphism $\cP^m(\R^p)\to \cP^m(\Gamma)$.
This justifies the notation $\cP^m$ for $\cP^m(\gG)$.

We have the obvious inclusion
\begin{equation}
    \cP^m\subset  \tilde\sym^m(\Gamma).
\end{equation}
Moreover, if $m<0$ then
\begin{equation}
     \sym^m(\Gamma)\cap \cP=\{0\}.
\end{equation}

If $\cA\subset\cinf{\Gamma}$ is a vector space which is closed under
$\pl_j,$  $j=1,...,p,$ we put
\begin{equation}
   \Omega^l\cA :=\Big\{\sum f_I dx_I\in\Omega^l(\Gamma)\,\Big|\, 
   f_I\in\cA\Big\}.
   \label{G1-2.1}
\end{equation}
If $\cA$ is graded then $\Omega^l\cA$ is bigraded, namely
\begin{equation}
  \Omega^l\cA^m:=\Big\{\sum f_I dx_I\in \Omega^l(\Gamma)\,\Big|\, 
f_I\in\cA^m\Big\}.
\end{equation}
Since $\cA$ is closed under $\pl_j$, $j=1,...,p$, the exterior derivative maps
$\Omega^l\cA$ into $\Omega^{l+1}\cA$, hence we obtain a complex
$(\Omega^*\cA,d)$. Obvious examples are
$\sym^*(\gG), \tilde\sym^*(\gG), \cP^*(\gG)$. 
Note that $d(\Omega^l\cP^m)\subset \Omega^{l+1}\cP^{m-1}$.

\begin{lemma} \label{S1-3.1} 
 The homology of the complex $(\Omega^*\cP,d)$ is given by
 $$  H^l(\Omega^*\cP,d)=\casetwo{\C}{l=0,}{0}{l\ge 1.}$$
 More precisely, if $\omega\in \Omega^l\cP^m$ with $l\ge 1$, is closed then
 there exists $\eta\in\Omega^{l-1}\cP^{m+1}$ with $d\eta=\omega$.
\end{lemma}
\begin{proof} $H^0(\Omega^*\cP,d)=\C$ is obvious. We mimick the usual proof
of the Poincar{\'e}\ Lemma (cf. \cite[Sec. I.4]{BotTu:DFAT}) and proceed
by induction on $p=\dim\Gamma$. If $p=1$ and $\omega=f(x)dx, f\in\cP^m$,
then we put $\eta(x):=\int_0^xf(t) dt\in\cP^{m+1}$. Hence the assertion
is true for $p=1$.

Next we consider
\begin{equation}
  \R^p\times \R \begin{array}{l}
  \ueber{\pi}{\longrightarrow}\\ \ueber{\longleftarrow}{s}
  \end{array} \R^p,\quad \pi(x,t)=x,\: s(x)=(x,0).
\end{equation}
Obviously,
\begin{equation}
\begin{split}
   & \pi^*(\Omega^l\cP^m(\R^p))\subset \Omega^l\cP^m(\R^p\times\R),\\
   & s^*(\Omega^l\cP^m(\R^p\times\R))\subset \Omega^l\cP^m(\R^p)
\end{split}
\end{equation}
holds. Since $\pi\circ s={\rm id}_{\R^p}$ we have $s^*\circ
\pi^*={\rm id}_{\Omega^l\cP^m(\R^p)}$.
An inspection of the construction of the usual homotopy operator $K$
(cf. \cite[Sec. I.4]{BotTu:DFAT}) shows that it induces a map 
\begin{equation}
   K:\Omega^l\cP^m(\R^p\times \R)\longrightarrow
   \Omega^{l-1}\cP^{m+1}(\R^p\times \R).
\end{equation}
Furthermore, $K$ satisfies the identity
\begin{equation}
   \id - \pi^*\circ s^* =dK+Kd
\end{equation}
from which  the assertion follows immediately.
\end{proof}

\begin{prop} \label{S1-3.2}Let $\go\in\Omega^1\tilde \sym^m(\Gamma), d\go=0$. 
If the cohomology class of $\omega$ in $H^1_{\rm dR}(\Gamma)$ vanishes
then there exists $F\in \tilde\sym^{m+1}(\Gamma)$ with $dF=\go$.
If $m+1<0$, then $F$ is uniquely determined.
\end{prop}
\begin{proof}
It is clear by assumption that $\go$ has a primitive in
$\cinf{\Gamma}$. The point is to prove that one of its primitives already
lies in $\tilde\sym^{m+1}(\Gamma)$. We write
$$\go=:\sum_{j=1}^p h_j dx_j.$$

If $m+1<0$, then the constant functions do not belong to
$\tilde\sym^{m+1}(\gG)$ which proves the uniqueness statement. 

If $m+1>0$, we fix $x_0\in\gG$ and put
\begin{equation}
  F(x):=\int_{x_0}^x\go,
\end{equation}
where $\int_{x_0}^x$ denotes integration along any path from $x_0$ to $x$.
This makes sense since $[\go]_{H^1_{\rm dR}(\gG)}=0$.
Since $\pl_jF=h_j\in\tilde\sym^m(\gG)$ it suffices to prove  the estimate
\begin{equation}
   |F(x)|\le C_{\eps,L} (1+|x|)^{1+\eps+m},\quad x\in L^c,
\end{equation}
for 
compact $L\subset \Gamma$ and 
$\eps>0$ (resp. $\eps=0$ if $m+1\not\in \Z_+$).

Since $h_j\in\tilde\sym^m(\gG)$ we have
$$|h_j(x)|\le C_{\eps,L} (1+|x|)^{m+\eps}$$
and thus
 $$  |F(x)|= \big| F(\frac{x}{|x|}) +\int_{\frac{x}{|x|}}^x \go\big|
     \le C+\big|\int_{\frac{x}{|x|}}^x \go\big|. 
 $$
Now
 $$  \int_{\frac{x}{|x|}}^x \go=\sum_{j=1}^p \int_{1/|x|}^1 h_j(tx) x_jdt
 $$
and consequently
\[ \begin{split}
   \big|\int_{\frac{x}{|x|}}^x \go\big|&\le
   C_{\eps,L}\int_{1/|x|}^1 (1+t|x|)^{m+\eps} |x|dt\\
   &\le  C_{\eps,L}(1+|x|)^{m+\eps+1}.\end{split}
\]
If $m=-1$, then $h_j(x)\le C(1+|x|)^{-1}$. Similarly
\begin{eqnarray*}
    |F(x)|&\le& C(1+ \int_{1/|x|}^1 (1+t|x|)^{-1} |x| dt) \le
    C_1+C_2\big|\log|x|\big|\\
    &\le& C_3 (1+|x|)^\eps
\end{eqnarray*}
holds for any $\eps>0$, hence we reach the conclusion in this case.

It remains to consider the case $m+1<0$. Then we put
\begin{equation}
F(x):= -\sum_{j=1}^p \int_1^\infty h_j(tx) x_j dt.
\end{equation}
Since $m+1<0$ we may differentiate under the integral.  Taking
$\pl_j h_l=\pl_l h_j$ into account we find $\pl_jF=h_j$. Moreover
one easily checks
$$ |F(x)|\le C_L (1+|x|)^{m+1},$$
which implies $F\in \tilde\sym^{m+1}(\gG)$.
\end{proof}

The partial derivatives $\pl_j$, $j=1,...,p$, are well--defined on the
quotient space $\tilde\sym^*(\Gamma)/\cP$,
hence we can form the complex $(\Omega^*(\tilde\sym^*(\Gamma)/\cP),d)$,
which is obviously isomorphic to the quotient complex
$(\Omega^*\tilde\sym^*(\Gamma),d)/(\Omega^*\cP,d)$.

\begin{prop} \label{S1-3.3}Let $\gG\subset \R^p$ be a connected cone.
For $\omega\in\Omega^1(\tilde \sym^m(\Gamma)/\cP),$ 
$d\omega=0,$ 
$[\go]_{H^1_{\rm dR}(\Gamma)}=0$ there exists $F\in\tilde\sym^{m+1}(\Gamma)/\cP$ with
$dF=\go$. $F$ is the unique element in 
$\tilde\sym^*(\Gamma)/\cP$ with $dF=\go$.
\end{prop}
\begin{remark} Here $[\go]_{H^1_{\rm dR}(\Gamma)}=0$ means that
$[\go_1]_{H^1_{\rm dR}(\Gamma)}=0$ for a closed representative $\go_1\in
\gO^1\tilde S^m(\gG)$ of $\go$. In view of Lemma \plref{S1-3.1} we
then have $[\go_2]_{H^1_{\rm dR}(\Gamma)}=0$
for any closed representative $\go_2$ of $\go$. 
\end{remark}
\begin{proof} We first prove that 
$\go=\sum\limits_{j=1}^p f_j dx_j$ has a closed
representative in $\gO^1\tilde\sym^m(\Gamma)$. Namely, pick representatives
$g_j\in\tilde\sym^m(\gG)$ of $f_j$ and put
$$\tilde\go=\sum_{j=1}^p g_jdx_j.$$
Since $d\go=0$ we have $d\tilde\go\in\gO^2\cP^{m-1}$. Since $d\tilde\go$
is closed, in view of Lemma \plref{S1-3.1} there exists $\eta\in\gO^1\cP^m$
with $d\eta=d\tilde\go$. Thus $\go_1:=\tilde\go-\eta$ is a closed
representative of $\go$. 

By Proposition \plref{S1-3.2} there exists $F_1\in\tilde\sym^{m+1}(\gG)$
with $dF_1=\go_1$. Then $F:=F_1\mod \,\cP$ sa\-tis\-fies
$dF=\go_1\mod\,\cP=\go$, which proves the existence of $F$.

If $\go_2\in\gO^1\tilde\sym^m(\gG)$ is another closed
representative of $\go$ and $F_2\in\tilde\sym^{m+1}(\gG), dF_2=\go_2,$
then 
%\begin{equation}
$$d(F_1-F_2)=\go_1-\go_2\in\gO^1\cP^m.$$
%\end{equation}
Since $d(\go_1-\go_2)=d^2(F_1-F_2)=0$ we again invoke Lemma \plref{S1-3.1}
and find a polynomial $P\in\cP^{m+1}$ with $dP=\go_1-\go_2$.
Hence $d(F_1-F_2-P)=0$ and thus $F_1-F_2=P+c\in\cP^{m+1}$.
This  proves the uniqueness statement.
\end{proof}

\begin{theorem} \label{S1-3.4} Let $M$ be a compact manifold, $E$ a smooth
vector bundle over $M$, and $\Gamma\subset \R^p$ a connected
cone with nonempty interior and
$H^1_{\rm dR}(\gG)=0$. Then there exists a unique linear map
$$\TR:{\rm L}^*(M,E,\gG)\rightarrow \sym^*(\gG)/\cP$$
with the following properties:
\begin{enumerate}\renewcommand{\labelenumi}{{\rm (\roman{enumi})}}
\item $\TR(AB)=\TR(BA)$, i.e. $\TR$ is a ``trace'',
\item $\TR(\pl_j A)=\pl_j \TR(A)$, \quad $j=1,...,p$,
\item If $A\in {\rm L}^m(M,E,\gG)$ and $m+\dim M<0$ then
$$\TR(A)(\mu)=\trltwo(A(\mu)).$$
\end{enumerate}

\noindent This unique $\TR$ satisfies furthermore:
\begin{enumerate}\renewcommand{\labelenumi}{{\rm (\roman{enumi})}}
\addtocounter{enumi}{3}
\item $\TR(\mu_j A)=\mu_j\TR(A)$, \quad $j=1,...,p$,
\item $\TR({\rm L}^m(M,E,\gG))\subset \tilde \sym^{m+\dim M}(\gG)$.
\end{enumerate}
\end{theorem}
\begin{remark} By slight abuse of notation $\mu_j$ denotes the operator
of multiplication by the $j$--th coordinate function. Note that
$\pl_j$ and $\mu_j$ is well--defined on the quotient
$\tilde\sym^*(\gG)/\cP$ since both operators map $\cP$ into itself.
Furthermore,
\begin{equation}
\begin{split}
   &\pl_j:{\rm L}^m(M,E,\gG)\rightarrow {\rm L}^{m-1}(M,E,\gG)\\
   &\mu_j:{\rm L}(M,E,\gG)\rightarrow  {\rm L}^{m+1}(M,E,\gG).
\end{split}
\end{equation}

\end{remark}
\begin{proof} 
\emph{Uniqueness:}\quad Assume there are $T_1,T_2$ satisfying
(i)--(iii). By (iii) $T_1$ and $T_2$ coincide on ${\rm L}^m(M,E,\gG)$ for
$m<-\dim M$. By induction assume that $T_1, T_2$ coincide on 
${\rm L}^{m_0}(M,E,\gG)$ and let $A\in {\rm L}^m(M,E,\gG), m\le m_0+1$, 
be given.

Consider the $1$--form
\begin{equation}
\begin{split}
   \go &:=\sum_{j=1}^p T_1(\pl_jA)(\mu) d\mu_j
       =\sum_{j=1}^p T_2(\pl_j A)(\mu) d\mu_j \,
       \in\gO^1(\tilde\sym^{m-1}(\gG)/\cP).
\end{split}
\end{equation}
In view of (ii) $\go$ is closed and we have
\begin{equation}
   dT_1(A)=\go = d T_2(A).
\end{equation}
Hence Proposition \plref{S1-3.3} implies $T_1(A)=T_2(A)$.

Next we assume that we have the unique $\TR$ with (i)--(iii) and prove that it
also satisfies (iv), (v).

If $m+1+\dim M<0$ and $A\in {\rm L}^m(M,E,\gG)$ then by (iii)
\begin{equation}
  \label{G1-3.5}
   \TR(\mu_j A)(\mu) =\trltwo(\mu_j A(\mu))=\mu_j \trltwo(A(\mu))=
    \mu_j \TR(A)(\mu)
\end{equation}
and
\begin{equation}
   \TR(A)=\trltwo(A(\cdot))\in \tilde\sym^{m+\dim M}(\gG)
   \label{G1-3.6}
\end{equation}
in view of Theorem \plref{S1-1.2}.

By induction we assume that \myref{G1-3.5} and \myref{G1-3.6} are true for 
$m\le m_0$. 
Now pick $A\in {\rm L}^m(M,E,\gG), m\le m_0+1$. Then
\begin{eqnarray}
    d \TR(\mu_j A)&=&\sum_{l=1}^p \TR(\mu_j \pl_l A) d\mu_l + \TR(A)
                     d\mu_j\nonumber\\
        &=&  \sum_{l=1}^p \mu_j \pl_l\TR(A) d\mu_l + \TR(A) d\mu_j\\
        &=& d(\mu_j \TR(A))\nonumber
\end{eqnarray}
and again by Proposition \plref{S1-3.3} we find $\TR(\mu_j A)=\mu_j \TR(A)$.

Similar, if $\TR(\pl_j A)\in \tilde \sym^{m-1+\dim M}(\Gamma)/\cP$ then 
\begin{equation}
  d\TR(A)=\sum_{j=1}^p \TR(\pl_j A) d\mu_j
\end{equation}
implies $\TR(A)\in\tilde\sym^{m+\dim M}(\gG)/\cP$. This proves (v) again by
induction.

\emph{Existence:}\quad Existence is also proved by induction
using Proposition \plref{S1-3.3}. Assume we have constructed $\TR$ on
${\rm L}^{m_0}(M,E,\gG)$. For $A\in {\rm L}^m(M,E,\gG), m\le m_0+1$, we let
$\TR(A)\in \tilde\sym^{m+\dim M}(\gG)/\cP$ be the unique primitive of the
closed $1$--form
\begin{equation}
   \go:=\sum_{j=1}^p \TR(\pl_j A)d\mu_j\in 
      \Omega^1(\tilde\sym^{m-1+\dim M}(\gG)/\cP).
\end{equation}
Obviously,  in this way we obtain a linear map $\TR:{\rm
L}^{m_0+1}(M,E,\gG)\rightarrow
\tilde \sym^*(\gG)/\cP$.

Now let $A\in {\rm L}^m(M,E,\gG), B\in {\rm L}^{m'}(M,E,\gG), m+m'\le m_0+1$.
Then
\begin{equation}\begin{split}
   &d\big[ \TR(AB)- \TR(BA)\big]\\
  =&\sum_{j=1}^p\big[ \TR((\pl_jA)B)+\TR(A\pl_jB)-
    \TR((\pl_j B)A)-\TR(B\pl_jA)\big] d\mu_j=0,    
                \end{split}
\end{equation}
hence $\TR(AB)=\TR(BA)$ again by Proposition \plref{S1-3.3}.

Since (ii) is obvious by construction we reach the conclusion.
\end{proof}
\begin{remark}
   Similarly to Theorem \plref{S1-1.3} $\TR$ is given by integration
of a canonically defined ``density''. See Remark \plref{ML5-4.4}
below.
\end{remark}

The construction of $\TR$ is much simpler, and more concrete, if $\Gamma$ is 
star-shaped:
\begin{proposition}
\label{S1-3.5}
 Let $\Gamma \subset \R^p$ be a star-shaped  cone with star-point $\mu_0$.
 Given $A \in {\rm L}^m (M,E,\Gamma)$ we have for fixed $\mu \in \Gamma$
 \begin{equation}
 \label{G1-3.7}
   A(\mu) - \sum_{|\alpha| \leq N-1} \, \frac{(\partial_\mu^\alpha
   A)(\mu_0)}{\alpha !} \, (\mu - \mu_0)^\alpha \in {\rm L}^{m-N} (M,E),
 \end{equation}
 and
 \begin{equation}
 \label{G1-3.8}
   \TR(A) \, (\mu) = \trltwo \Big( A(\mu) - \sum_{|\alpha| \leq N-1}
   \, \frac{(\partial_\mu^\alpha A) (\mu_0)}{\alpha !} \, (\mu - \mu_0)^\alpha
   \Big) \, {\rm mod} \, {\cP}  ,
 \end{equation}
 where $N$ is large enough.
\end{proposition}
\begin{proof}
  Taylor's formula implies
  \[\begin{split}
  \lefteqn{ A(\mu) - \sum_{|\alpha| \leq N-1} \, \frac{(\partial_\mu^\alpha
    A)(\mu_0)}{\alpha !} \, (\mu - \mu_0)^\alpha}\\
    & =\frac{1}{(N-1)!} \, \sum_{|\alpha| = N} \, 
    \int_0^1 \, (1 -t)^{N-1} \,
    (\partial_{\mu}^{\alpha} A) (\mu_0 + t (\mu -\mu_0) ) \, dt \;
     (\mu -\mu_0)^\alpha \in {\rm L}^{m-N} (M,E),
  \end{split}
  \]%\end{equation}
  hence the right hand side of \myref{G1-3.8} is well-defined. 
  Now it is easy to check that the right hand side of \myref{G1-3.8} 
  defines a (well--defined, i.e. independent of $N$) map with the 
  properties (i), (ii) and (iii) of Theorem \ref{S1-3.4}. By uniqueness
  we reach the conclusion.
\end{proof}

We tried hard to prove that for the algebra
${\rm CL}(M,\Gamma)$ the 
properties (i), (ii), and (iv)
of Theorem \plref{S1-3.4} already
determine $\TR$ up to a scalar factor. This would be nicer than 
assuming (iii),
which prescribes $\TR$ on a large class of operators.

We state this as a conjecture:

\begin{conjecture}\label{ML5-3.9}
Let $M$ be a compact manifold and $\Gamma\subset\R^p$ be a connected
cone with nonempty interior and
$H^1_{\rm dR}(\Gamma)=0$. Let
$\tau:{\rm CL}^*(M,\Gamma)\rightarrow \tilde\sym^*(\Gamma)/\cP$
be a linear map satisfying {\rm (i), (ii),} and {\rm (iv)} of
Theorem \plref{S1-3.4}.
Then there is a constant $c$ such that $\tau=c \,\TR$.
\end{conjecture}
But we have a partial result:
\begin{prop}\label{ML5-3.10} Let $M$ and $\Gamma$ be as in the 
preceding conjecture. Let $\tau:{\rm L}^{-\infty}(M,\Gamma)
\rightarrow \sym^{-\infty}(\Gamma)$
be a linear map satisfying   {\rm (i), (ii),} and {\rm (iv)} of
Theorem \plref{S1-3.4}.
Then there is a constant $c$ such that 
$\tau=c\, \TR\restr {\rm L}^{-\infty}(M,\Gamma)$, 
i.e. for 
$A\in {\rm L}^{-\infty}(M,\Gamma)$ 
$$\tau(A)(\mu)=c \,\trltwo(A(\mu)).$$
\end{prop}
\begin{proof}
$\tau$ is local, e.g. $\tau(A(\cdot))(\mu)=0$ if $A(\mu)=0.$
To see this let $A(\mu_0)=0$. We write
$$A(\mu)= \sum_{j=1}^p A_j(\mu)(\mu_j-\mu_{0,j})$$
with $A_j\in {\rm L}^{-\infty}(M,\Gamma).$ Then in view of (iv)
we have 
$\tau(A)(\mu_0)=0$.

Next we pick $\mu_0\in\Gamma$ and choose a function
$f\in \cinfz{\Gamma}$ with $f(\mu_0)=1$.
For $K\in {\rm L}^{-\infty}(M)$ we have
$f(\cdot)K\in {\rm L}^{-\infty}(M,\Gamma)$. Consequently,
\begin{equation}
     \tau_{\mu_0}:{\rm L}^{-\infty}(M)\rightarrow \C, \quad
   K\mapsto \tau(f(\cdot)K)(\mu_0)\label{ML5-3.24}
\end{equation}
is a trace on ${\rm L}^{-\infty}(M)$.

Indeed, from the locality of $\tau$ we conclude that $\tau_{\mu_0}$ is
independent of the $f$ chosen and in view of (i) $\tau_{\mu_0}$
is a trace.

Since each trace on ${\rm L}^{-\infty}(M)$ is a scalar multiple
of the $L^2$--trace (see e.g. \cite[Appendix]{Gui:RTCAFIO})
there is a constant $c(\mu_0)$ such that
\begin{equation}
      \tau_{\mu_0}=c(\mu_0) \trltwo.
     \label{ML5-3.25}
\end{equation}

For $f\in\cinfz{\Gamma}$ with $f(\mu_0)\neq 0$ we find in view
of \myref{ML5-3.25}
\begin{equation}
     \tau(f(\cdot)K)(\mu_0)= f(\mu_0)c(\mu_0) \trltwo(K).
   \label{ML5-3.26}
\end{equation}
By the locality of $\tau$
\myref{ML5-3.26} also holds if $f(\mu_0)=0$.
Since $\tau(f(\cdot)K)\in\sym^{-\infty}(\Gamma)$ we infer from
\myref{ML5-3.26} that the function $\mu\mapsto c(\mu)$ is smooth.

Now loc. cit. gives for $f\in\cinfz{\Gamma}$
\begin{equation}
    \pl_j \tau(f(\cdot)K)=((\pl_j f) c+f\pl_j c)\trltwo(K).
\end{equation}
On the other hand using (ii) and \myref{ML5-3.26}
\begin{equation}
    \pl_j \tau(f(\cdot)K)=\tau((\pl_jf)(\cdot)K)=c (\pl_jf) \trltwo(K).
\end{equation}
Hence $\pl_j c=0$, $j=1,...,p$, and thus $c$ is a constant.

Finally we invoke again the locality of $\tau$. Let 
$A\in{\rm L}^{-\infty}(M,\Gamma)$
and $f\in\cinfz{\Gamma}$ with $f(\mu_0)=1$. Then
\begin{equation}
   \begin{split}
     \tau(A)(\mu_0)=& \tau(A- f(\cdot) A(\mu_0))(\mu_0)+\tau(f(\cdot)A(\mu_0))(\mu_0)\\
  =& \tau_{\mu_0}(A(\mu_0))=c\,\trltwo(A(\mu_0)),
   \end{split}
\end{equation}
since $A- f(\cdot) A(\mu_0)$ vanishes at $\mu_0$.
\end{proof}

This proof cannot be extended to prove the Conjecture \plref{ML5-3.9},
at least not in an obvious way. The main reason is that
if $K\in 
{\rm L}^m(M,\Gamma)$ (resp. ${\rm CL}^m(M,\Gamma)$) , 
$m>-\infty,$ and $f\in\cinfz{\Gamma}\setminus\{0\}$
then $f(\cdot)K\not\in {\rm L}^*(M,\Gamma)$ (resp. ${\rm CL}^m(M,\Gamma)$).

Another idea of proving Conjecture \plref{ML5-3.9} is to mimick 
the uniqueness proof for the noncommutative residue.
If every $A\in{\rm CL}^m(M,\Gamma)$ could be written 
\begin{equation}
  A=\sum_{j=1}^N [P_j,Q_j]+R,
  \label{ML5-3.30}
\end{equation}
with $R\in {\rm L}^{-\infty}(M,\Gamma)$ and $P_j,Q_j\in {\rm CL}(M,\Gamma)$
(as it is the case for $\Gamma=\{0\}$) then the Conjecture \plref{ML5-3.9}
would immediately follow from the previous Proposition \plref{ML5-3.10}.

However, there is some evidence that \myref{ML5-3.30} actually is wrong
in general.

%%% Local Variables: 
%%% mode: latex
%%% TeX-master: "LesPflTAPDPOEI"
%%% End: 

%
\section{Exotic traces on $\rm\bf CL^*(M,E,\R^p)$}
Let $\TR$ be the map defined in Theorem \plref{S1-3.4}. 
First we want to characterize the space $\TR(\CL^*(M,E,\gG))$.
We know
from Theorem \plref{S1-1.2} that $\TR(\CL^m(M,E,\gG))\subset \CS^{m+\dim M}(\gG)$
if $m+\dim M<0.$

\begin{lemma} \label{S1-4.1} Let $\gG\subset \R^p$ be a connected cone
with nonempty interior and
$H^1_{\rm dR}(\gG)=0$. If $A\in\CL^m(M,E,\gG)$ then
$\TR(A)$ has a representative $f\in \sym^*(\gG)$ which has an asymptotic
expansion in $\sym^*(\gG)$
\begin{equation}
\label{G3-4.1}
  f\sim \sum_{j=0}^\infty f_j + \sum_{k=0}^{[m+\dim M]} \, g_k
  +\sum_{k=0}^{[m+\dim M]}  h_k,
\end{equation}
where $f_j(\mu)$ is homogeneous of degree $m-j$ for $|\mu|\ge 1,$
$g_k(\mu)$ is homogeneous of degree $k$ for $|\mu|\ge 1$, and
$h_k(\mu)=\tilde{h}_k(\mu/|\mu|)\,|\mu|^k \log|\mu|$ for $|\mu|\ge 1$
and a smooth function 
$\tilde{h}_k \in C^\infty (\Gamma \cap S^{p-1})$.
\end{lemma}
\begin{proof}
  If $m + {\rm dim} \, M < 0$ then the assertion follows from Theorem 
  \ref{S1-1.2}. Assume by induction that the assertion is true for 
  $m \leq m_0$ and consider $A \in {\rm CL}^m (M,E, \Gamma)$ with 
  $m \leq m_0 +1$. Then $\TR (\partial_j A)$, $j =1,...,p$, has a 
  representative of the form \myref{G3-4.1}. Integrating the 
  representative of $d \, \TR (A)$ we reach the conclusion.
\end{proof}

\begin{dfn} \label{S2-4.2} Let $\gG\subset \R^p$ be a cone. We denote by
${\rm PS}^m(\gG)$ the set of all functions $a\in \sym^*(\Gamma)$ 
which admit an asymptotic expansion
 $a\sim \sum\limits_{j\ge 0} a_{m_j}$ in $\sym^*(\gG)$, where
$a_{m_j}\in\cinf{\Gamma}$, $m\ge m_j\searrow -\infty$ and
$$a_{m_j}=\sum_{l=0}^{k_j} g_{lj},$$
with $g_{lj}(\xi)=\tilde{g}_{lj}(\xi/|\xi|)|\xi|^{m_j} \log^l|\xi|$
for $|\xi|\ge 1$.
We call these functions \emph{$\log$--poly\-homo\-geneous}
(cf. \cite{Les:NRPOPS}).
\end{dfn}
By Lemma \plref{S1-4.1} $\TR(A)$ is $\log$--polyhomogeneous 
for $A\in\CL^*(M,E,\gG)$.

From now on we content ourselves to $\gG=\R^p$. First we are going to
introduce a regularized integral for $\log$--polyhomogeneous 
functions. Let us consider a $\log$--polyhomogeneous function
$f$ and write
$$f=\sum_{m_j\ge-N} f_{m_j}+g,$$
where
$$g(\xi)=O(|\xi|^{-N}),\quad |\xi|\to\infty.$$
Thus we have an asymptotic expansion
\begin{equation}
    \int_{|\xi|\le R} f(\xi)d\xi\sim_{R\to\infty}
       \sum_{\ga\to-\infty}
       p_\ga(\log R) R^\ga, \label{G2-4.1}
\end{equation}
where $p_\ga$ is a polynomial of degree $k(\ga)$. 
Then we define $\reginttext_{\R^p} f(x) dx$ to be the constant term in
this asymptotic expansion, i.e.
\begin{equation}
\regint_{\R^p} f(x)dx:= \LIM_{R\to \infty} \int_{|x|\le R} f(x) dx
   :=p_0(0),
   \label{ML5-4.3}
\end{equation}
(cf. \cite[Sec. 5]{Les:NRPOPS}). 
If $p=1$ then $\reginttext$ coincides with the Hadamard partie finie.

Here and in the sequel we will
use the common notation $\LIM\limits_{R\to\infty}$ for the
constant term in the $\log$--polyhomogeneous expansion as $R\to
\infty$.

The transformation behavior of $\reginttext_{\R^p}$ is more complicated
than that of the usual integral.
The following result gives the change of variables formula for $\reginttext$
under linear coordinate changes:

\begin{prop}\label{S2-4.3}\hspace*{-2pt}{\rm (cf. \cite[Prop.~5.2]{Les:NRPOPS})}
Let $A\in {\rm GL}(n,\R)$ be a regular matrix. Further\-more, let
$f \in{\rm PS}^*(\R^p)$ with
$f
\sim \sum\limits_{\ga\to-\infty} f_\ga$,
$f_{\ga}(\xi)=:\sum\limits_{l=0}^{k(\ga)} f_{\ga,l},
\;  f_{\ga,l}(\xi)=f_{\ga,l}(\xi/|\xi|)|\xi|^{\ga}
\log^l(|\xi|)$, $|\xi|\ge 1$. Then we have
$$\regint_{\R^p} f(A\xi) d\xi=
   |\det A|^{-1}\left( \regint_{\R^p} f(\xi)d\xi+
      \sum_{l=0}^{k}\frac{(-1)^{l+1}}{l+1}\int_{S^{n-1}} f_{-n,l}(\xi)
          \log^{l+1} |A^{-1}\xi| d\xi\right).$$
\end{prop}
\begin{remark}\label{ML5-4.4}
Proposition \plref{S2-4.3} allows to give an alternative existence proof
for $\TR$ which in addition shows that $\TR$ is given by integration
of a canonical density.
Furthermore, this alternative construction reveals the analogy to the
Kontsevich-Vishik trace.

We consider the local situation. For simplicity let $E=\C$, $U$ a 
coordinate patch and $A\in\CL^m(U,\Gamma)$ as in \myref{ML5-1.3}.
For fixed $x\in U,\mu\in\Gamma$ the symbol $\xi\mapsto
\sigma_A(x,\xi,\mu)$ is polyhomogeneous. We define
the density
\begin{equation}
    \go_A(x,\mu):=\regint_{\R^n} \sigma_A(x,\xi,\mu)\dbar\xi |dx|
      \mod \cinfty{U}\otimes\cP.
\end{equation}
Using Proposition \plref{S2-4.3} one shows similarly to
\cite[Lemma 5.3]{Les:NRPOPS} that in this way one obtains a well--defined
``density''
\begin{equation}
      \go_A\in 
    \tilde\sym^{m+\dim M}(M,|\Omega|,\Gamma)\Big/\cinfty{M,|\Omega|}\otimes \cP
\end{equation}
fulfilling
\begin{equation}
      \TR(A)=\int_M\go_A.
\end{equation}
Note that $\int_M\go_A$ is a well--defined element of 
$\tilde\sym^{m+\dim M}(M,\Gamma)/\cP$.
Thus we obtain the analogue of Theorem \plref{S1-1.3} for arbitrary $m$.
The details are left to the reader.

In \cite[Sec.~5]{Les:NRPOPS} the first named author used this kind of argument
to construct the Kontsevich-Vishik trace.
\end{remark}

Let us mention that Stokes' theorem does not hold for $\reginttext$, or
in other words $\reginttext$
is not a closed functional on $\gO^*({\rm PS}(\R^p))$. 
More precisely, we extend $\reginttext$ to $\gO^* ({\rm PS}^*(\R^p))$
by putting 
\begin{equation}
  \regint:\go\mapsto\casetwo{0}{\go\in \gO^k,k<p,}{
    \regint_{\R^p}f(\xi)d\xi}{\go=f(\xi) d\xi_1\wedge\ldots\wedge d\xi_p.}
\end{equation}
In this way we obtain a graded trace on the complex 
$(\gO^* ({\rm PS}^*(\R^p)),d)$. This would be a cycle in the sense of
{\rm Connes}  \cite[Sec. III.1.$\alpha$]{Con:NG} if $\reginttext$
were closed. 

The next lemma shows that $d\reginttext$, which
is defined by $(d\reginttext)\go:=\reginttext d\go$, 
is nontrivial. However, it is local in the sense that it depends
only on the $\log$--polyhomogeneous expansion of $\go$.

\begin{lemma} \label{S2-4.4} Let $f\in {\rm PS}^*(\R^p)$, 
$f\sim \sum\limits_{\ga\to -\infty}f_\ga$, where 
$$f_\ga(\xi)=\sum_{l=0}^{k(\ga)}f_{\ga,l}(\xi/|\xi|) |\xi|^\ga
\log^l |\xi|, \quad |\xi|\ge 1, \quad 
f_{\ga,l}\in\cinf{S^{p-1}}.$$ Then
$$\regint_{\R^p} \frac{\pl f}{\pl \xi_j} d\xi=
   \int_{S^{p-1}} f_{1-p,0}(\xi) \xi_j d{\rm vol}_S(\xi).$$
\end{lemma}
\begin{proof}
It suffices to prove this formula for $f\in \cinf{\R^p}$ with
$$f(\xi)= f(\xi/|\xi|) |\xi|^\ga \log^l|\xi|,\quad |\xi|\ge 1.$$
Then by Gau{\ss}' formula
\begin{equation}\begin{split}
     \int_{|\xi|\le R}\frac{\pl f}{\pl \xi_j}\, d\xi&=
      \int_{|\xi|=R} f(\xi) \frac{\xi_j}{|\xi|}\, d{\rm vol}_S(\xi)\\
      &= \int_{|\xi|=1} f(\xi) \xi_j\, d{\rm vol}_S(\xi)\:
            R^{p-1+\ga} \log^l R
                \end{split}
\end{equation}
and we reach the conclusion.
\end{proof}

Next we consider $f\in \cP$. As
$f$ is a sum of homogeneous polynomials
\begin{equation}
   f=\sum_{k=0}^m f_k,\quad f_k(\gl\xi)=\gl^k f_k(\xi)
\end{equation}
we have $\cP\subset {\rm PS}^*(\R^p)$.
Thus
\begin{equation}
    \int_{|\xi|\le R} f(\xi) d\xi=\sum_{k=0}^m \int_{|\xi|\le 1}
  f_k(\xi) d\xi\: R^{k+p}
\end{equation}
and hence
\begin{equation}
       \regint_{\R^p} f(\xi) d\xi =0.
\end{equation}
As a consequence $\reginttext$ factorizes through $\cP$ to
a well--defined functional on ${\rm PS}^*(\R^p)/\cP$.

%Using $\reginttext$ we can define a trace on $\CL^*(M,E,\R^p)$:
%
%Since $\regint_{\R^p}$ vanishes on polynomials, $\overline{\Tr}$ is
%well--defined and is a trace on $\CL^*(M,E,\gG)$. We leave it to the 
%reader to check that in the case $p=1$ this extended traces coincides
%with the extended trace defined in \cite{Mel:EIFPO}.
%
%Next we want to explain how the formal trace comes into play:

For an arbitrary connected cone $\Gamma$ with $H^1_{\rm dR}(\gG)=0$
we can construct a complex from
$\CL^*(M,E,\gG)$. Namely, similarly to
\myref{G1-2.1} we put
\begin{equation}
   \gO^k\CL^*(M,E,\gG):=\Big\{\sum A_I dx_I\,\Big|\, |I|=k,
   A_I\in \CL^*(M,E,\gG)\Big\}.
   \label{pseudo-G1-4.1}
\end{equation}

The exterior derivative maps  $\gO^*\CL^*(M,E,\gG)$ into itself and
so we obtain a complex $(\gO^*\CL^*(M,E,\R^p),d)$. The cup product makes
$\gO^*\CL^*(M,E,\R^p)$ into a graded algebra and $\TR$ extends naturally
to a complex homomorphism
\newcommand{\Trbar}{\overline{\Tr}}
\begin{equation}\begin{array}{rrll}
  \TR:&(\gO^*\CL^*(M,E,\gG),d)&\longrightarrow &(\gO^*{\rm PS}^*(\gG)/\cP,d)\\[1em]
   &\sum A_I dx_I&\longmapsto &\sum \TR(A_I) dx_I.
                \end{array}
\end{equation}
That $\TR$ is indeed a complex homomorphism follows from
Theorem \plref{S1-3.4}.

If $\tau:{\rm PS}^*(\gG)\to \C$ is any linear functional with
$\tau|\cP=0$ then $\tau$ factorizes through $\cP$ to a linear functional
on ${\rm PS}^*/\cP$ and we obtain a graded trace on 
$(\gO^*{\rm PS}^*(M,E,\gG)/\cP,d)$ by putting
\begin{equation}
    \bar \tau(\go):=\casetwo{0}{\go\in \gO^k,k<p,}{
    \tau(f)}{\go=fdx_1\wedge\ldots\wedge dx_p.}
\end{equation}
$\tau$ induces a graded trace on $(\gO^*{\rm CL}^*(M,E,\gG),d)$ by
putting $\bar\tau:=\tau\circ \TR$.
In particular we can apply this construction to $\reginttext$ and
$d\reginttext$:
\begin{dfn} \label{S1-4.2} On $\gO^*\CL^*(M,E,\R^p)$ we define the extended
trace
$$\overline{\Tr}:=\regint\circ \TR,$$%\regint_{\R^p} \TR(A)(\mu) d\mu,$$
and the formal trace
$$\widetilde\Tr:=d\Trbar:=\left( d\regint \right) \circ \TR=
  \regint\circ d\circ\TR =\regint\circ\TR\circ d.$$
\end{dfn}

\begin{remark}

\mbox{}
\begin{enumerate}
  \renewcommand{\labelenumi}{{\rm \arabic{enumi}.}}
\item
  The extended trace is graded, 
  thus $((\gO^*\CL^*(M,E,\R^p),d),\Trbar)$ is almost
  a cycle in the sense of {\sc Connes} \cite[Sec. III.1.$\alpha$]{Con:NG},
  except that $\Trbar$ is not closed. Its derivative, the formal trace, 
  $\widetilde \Tr$,
  is a closed graded trace on $\gO^*\CL^*(M,E,\R^p)$.
  Furthermore, Lemma \plref{S2-4.4} shows that it is symbolic, 
  i.e. depends only on finitely many terms of the symbol expansion of $A$
  (see the next proposition below). 
\item
  If $p =1$ then via the isomorphism of Proposition 
  \plref{ML1-S1.1} the traces
  $\overline{\Tr}, \frac{1}{2\pi}\widetilde{\Tr}$ 
  coincide with the 
  traces  $\overline{\Tr}, \widetilde{\Tr}$ 
  introduced by {\sc R.~B.~Melrose} 
  \cite[Sec.~4 and 7]{Mel:EIFPO}.
  Note that our normalization of the formal
  trace $\widetilde\Tr$ differs from the one of \cite{Mel:EIFPO}
  by a factor of $\frac{1}{2\pi}.$
 \end{enumerate}
\end{remark}
\begin{proposition}
  \label{S3-4.6} {\rm (cf.~\cite[Prop.~6]{Mel:EIFPO})}
  Let $A \in {\rm CL}^m (M, E, \R^p)$ and put
  \begin{equation}
    \omega := (-1)^{j-1} A \, d\mu_1 \wedge ... \wedge \widehat{d \mu_j} \wedge
    ... \wedge d \mu_p.
  \end{equation}
  Then 
  \begin{equation}
    \widetilde{\Tr} (\omega) = \lim_{L \rightarrow \infty} \int_{|\xi| \leq L}
    \int_{S^{p-1}} \, \tr (a_{1-p-n} (x,\xi,\mu) )\, \mu_j \, d{\rm vol} (\mu)
    \, dx \, \dbar \xi.
  \end{equation}
\end{proposition}
\begin{proof}
  It suffices to prove this for $E=\C$. Obviously, we have 
  \begin{equation}
    \widetilde{\Tr} (\omega) = \overline{\Tr} (d \omega) = 
    \regint_{\R^p} \TR \left( \partial_{\mu_j} A \right) \, d \mu .
  \end{equation}
  In view of Proposition \ref{S1-3.5} we introduce the abbreviation
  \begin{equation}
    b_N (x,\xi,\mu) := \sum_{|\alpha| \leq N-1} \frac{(\partial_\mu^\alpha)
    a (x,\xi,0)}{\alpha !} \, \mu^\alpha.
  \end{equation}
  Using Lemma \ref{S2-4.4} we find for $N$ large enough
%  \begin{equation}
 \[ \begin{split}
    \regint_{\R^p} \TR \left( \partial_{\mu_j}A \right) \, (\mu) \, d\mu & = 
    \LIM_{R \rightarrow \infty} \, \int_{|\mu|\leq R} \, \int_U \, \int_{\R^n}
    \left[ (\partial_{\mu_j} a)(x,\xi,\mu) - (\partial_{\mu_j} b_N )(x,\xi,\mu)
    \right] \dbar \xi \, dx \, d\mu \\
    & = \int_U \, \int_{S^{p-1}} \left\{ \int_{\R^n} \, a(x,\xi,\mu) -
    b_N(x,\xi,\mu) \, \dbar \xi  \right\}_{1-p} \mu_j \, d{\rm vol}_S (\mu)
    \, dx \\
    & =\lim_{L \rightarrow \infty} \, \int_{|\xi| \leq L} \int_{S^{p-1}}
    a_{1-p-n} (x,\xi,\mu)\, \mu_j \, d{\rm vol}_S (\mu) \, \dbar\xi \, dx.    
  \end{split}  
 \]% \end{equation}
  Here $\{ \, .\,\}_{p-1}$ denotes the term of $\mu$-homogeneity $1-p$.
\end{proof}

%%% Local Variables: 
%%% mode: latex
%%% TeX-master: "LesPflTAPDPOEI"
%%% End: 

%
\section{The eta--invariant}
Now we have all prerequisites to define the higher eta--invariants:
\begin{dfn}\label{S5-5.1}
  If $p = 2k-1$ then we put for elliptic and invertible 
  $A \in {\rm CL}^* (M, E,\R^{2k-1} )$
  \begin{equation}
    \eta_k (A) := 2 c_k \overline{\Tr} 
    \left( \left( A^{-1} dA\right)^{2k-1} \right),
  \end{equation}
  where $c_k = \frac{(-1)^{k-1} (k-1)!}{(2 \pi i)^k (2k-1)!}$.
\end{dfn}
\begin{remark}

\mbox{}
\begin{enumerate}\renewcommand{\labelenumi}{{\rm \arabic{enumi}.}}
\item  If $k=1$ then via the isomorphism of Proposition \plref{ML1-S1.1}
$\eta_1 (A) = \frac{1}{\pi i} \overline{\Tr} \left( A^{-1} dA \right)$
  coincides with the $\eta$-invariant of {\sc Melrose} 
  \cite[Sec.~5]{Mel:EIFPO}.
\item There are at least two
motivations for the choice of the normalization
constant $c_k$:
It is well--known that 
for every smooth map $f:S^{2k-1}\to {\rm GL}(N,\C)$ the number
$$w(f):=c_k\int_{S^{2k-1}}\tr((f^{-1}df)^{2k-1})$$
actually is an integer and $w$ induces an isomorphism 
$\pi_{2k-1}({\rm GL}(\infty,\C))\to\Z$.
A map with $w(f)=1$ can be constructed using Clifford matrices
(cf. \myref{ML5-5.23} below).
In this sense $\eta_k$ is a higher ``winding number''.

The second motivation comes from the relation to the 
spectral eta--invariant (see Proposition \plref{ML5-5.6} below).
\end{enumerate}
\end{remark}
\begin{proposition} {\rm (cf. \cite[Prop. 7]{Mel:EIFPO})}
  \label{ML5-5.3}  
Let $A_s \in {\rm CL}^* (M,E, \R^{2k-1} ) $ be elliptic invertible and
  smoothly dependent on $s \in [0,1]$. Then
  \begin{equation}
    \frac{d}{ds} \eta_k (A_s) = 2(2k-1) c_k \widetilde{\Tr} \left(
    (A_s^{-1} \partial_s A_s) (A_s^{-1} d A_s)^{2k-2} \right).
  \end{equation}
\end{proposition}
\begin{proof}
We introduce the $1$--form $\go:=A^{-1}dA$ and note that
since 
\begin{equation}
      d\go=-\go\wedge\go
\end{equation}
we have for $l\in \Z_+$
\begin{equation}
         d\go^l=\smallcasetwo{0}{l\;{\rm even},}{-\go^{l+1}}{l\;{\rm odd}.}
\end{equation}
Using this we find
  \begin{equation}
    \begin{split}
      \frac{d}{ds} \eta_k (A_s) & = 2c_k \, \overline{\Tr} 
      \left( \frac{d}{ds} (A_s^{-1} dA_s )^{2k-1} \right) \\
      & =2 (2k-1) \, c_k \, \overline{\Tr} \left( (A_s^{-1} dA_s )^{2k-2} 
      \frac{d}{ds} A_s^{-1} dA_s \right)\\
      & =2 (2k-1) \, c_k \, \overline{\Tr} \left( - (A^{-1} dA )^{2k-1}
      \, A_s^{-1} \partial_s A +  (A_s^{-1} dA_s )^{2k-2}  
        \, A_s^{-1} \, d \partial_s A_s \right) \\
      & = 2(2k-1) \, c_k \, \overline{\Tr} 
      \left( d \left[ (A_s^{-1} dA_s )^{2k-2}  \, A_s^{-1} \, \partial_s A_s 
      \right]\right) \\
      & = 2(2k-1) \, c_k \, \widetilde{\Tr} 
      \left( (A_s^{-1} \partial_s A_s ) (A_s^{-1} dA_s )^{2k-2} \right).
    \end{split}
  \end{equation}
\end{proof}
Next we consider the complex Clifford algebra 
$\C\ell_p$, $p=2k-1$ odd.
$\C\ell_p$ is the universal $C^*$--algebra generated by $p$ unitary elements
$e_1,...,e_p$ subject to the relations
\begin{equation}
    e_i \cdot e_j + e_j \cdot e_i = - 2 \delta_{ij}
    \label{Cliffrel}
\end{equation}
(cf. e.g. \cite[Chap. I]{LawMic:SG}).
We choose a Clifford representation 
\begin{equation}
  c: \R^p \rightarrow {\rm M}(N,\C), \quad c(x) = \sum_{j=1}^p \, x_j E_j,
\end{equation}
where the $E_j$ are skew-adjoint Clifford matrices in $\C^N$. This means
that $E_1,\ldots,E_p$ are skew--adjoint matrices satisfying the Clifford
relations \myref{Cliffrel}. $c$ induces a $*$--representation of
$\C\ell_p$ in ${\rm M}(N,\C)$.

Let us introduce the map 
\begin{equation}
              f:\R \times \R^p \rightarrow {\rm M} (N,\C),\; 
         x=(x_0,x') \mapsto x_0  + c(x')
   \label{G5-5.5}
\end{equation}
and the $p$-form
\begin{equation}
  \omega = \tr \left( f^{-1} \, df \right)^p
\end{equation}
on $\R^{p+1} \setminus \{ 0 \}$. That $f(x)$ is indeed invertible for
$x\neq 0$ follows from \myref{ML5-5.11} below.
\begin{proposition}
\label{MPS1-4.1}
  The $p$-form $\omega$ is given by
  \begin{equation}
    \omega = |x|^{-p-1} \, p! \, \tr (E_1 \cdot ... \cdot E_p)
    \sum_{j=0}^p \, (-1)^j \, x_j \, dx_0 \wedge ... \wedge \widehat{dx_j} \wedge
    ... \wedge dx_p.
  \end{equation}
\end{proposition}
\begin{proof}
  By definition of $f$ we have $f(x)^* f(x) = |x|^2 $, hence 
  \begin{equation}
    f(x)^{-1} = |x|^{-2} f(x)^*,
    \label{ML5-5.11}
  \end{equation}
  and
  \begin{equation}
    f(x)^* \, df (x) = (x_0 - c(x')) \, \left( dx_0 + \sum_{j=1}^p \, E_j 
    \, dx_j  \right).
  \end{equation}
  Setting $\omega_1 = (x_0 - c(x')) dx_0 $ and 
  $\omega_2 = (x_0 - c(x')) \sum\limits_{j=1}^p \, E_j dx_j $ we get  
  $\omega_1^2 =0$ and
  \begin{equation}
    \left( f(x)^* \, df(x)  \right)^p = \omega_2^p + \sum_{j=0}^{p-1} \,
    \omega_2^j \, \omega_1 \, \omega_2^{p-1-j}  ,  
  \end{equation}
  which implies 
  \begin{equation}
    \omega = |x|^{-2p} \, \left( \tr \omega_2^p  + 
    p \, \tr \omega_1 \, \omega_2^{p-1} \right) ,
    \label{G5-5.67}
  \end{equation}
  as $p$ is odd.
  We calculate the two terms on the right hand side separately. 
  First note that both terms are invariant with respect to 
  transformations of the form $x' \mapsto O x'$ with $O \in SO (p)$, 
  so we may assume $x' =(x_1,0,...,0)$. Now
  \begin{equation}
  \begin{split}
     &\tr \omega_2^p \big|_{(x_0,x')} = \tr \Big(
    ( x_0 -x_1 E_1) 
    \sum_{j=1}^p \, E_j \, dx_j \Big)^p \Big|_{(x_0,x')} \\
    =& \sum_{\sigma \in {\rm S}_p} \, ( \sgn \sigma ) \:
    \tr \left( (x_0 -x_1E_1) E_{\sigma (1)} \cdot ... \cdot 
    (x_0 -x_1E_1) E_{\sigma (p)} \right) \, dx_1 \wedge ... \wedge dx_p.
  \end{split}
  \end{equation}
  We fix a permutation 
  $\sigma \in {\rm S}_p$ and put $j:=\sigma^{-1} (1)$. Then for
  $l\neq j,j-1$ we have the relation
  \begin{equation}
    (x_0 -x_1 E_1) E_{\sigma (l)} (x_0 - x_1 E_1) E_{\sigma (l+1)}
    = | x |^2 E_{\sigma (l)} \, E_{\sigma (l+1)}.
  \end{equation}
  Hence we obtain
  \begin{equation}
    \begin{split}
    \lefteqn{\tr \left( (x_0 - x_1 E_1 ) E_{\sigma (1)} \cdot ... \cdot 
      (x_0 -x_1E_1) E_{\sigma (p)} \right) =}\\
      & = \tr \left( (x_0 - x_1 E_1 ) E_1 \, (x_0 - x_1 E_1 ) E_{\sigma (j+1)}
      \cdot ... \cdot (x_0 - x_1 E_1 ) E_{\sigma (p)} \cdot \right. \\
      & \hspace{25mm} 
      \left. \cdot (x_0 - x_1 E_1 ) E_{\sigma (1)} \cdot ... \cdot 
      (x_0 - x_1 E_1 ) E_{\sigma (j-1)} \right) \\
      & = |x|^{2(k-1)} \tr \,  \left( (x_0 - x_1 E_{\sigma (j)}) \cdot
      E_{\sigma (j+1)} \cdot ... \cdot E_{\sigma (p)} \cdot E_{\sigma (1)}
      \cdot ... \cdot E_{\sigma (j)} \right) \\
      & = x_0 |x|^{p-1} \, \tr \left( E_{\sigma (1)} \cdot ... \cdot
      E_{\sigma (p)} \right),  
    \end{split}
  \end{equation}
  where we have used the Berezin Lemma (cf.~\cite[Prop.~3.21]{BerGetVer:HKDO}).
  Summing up we have proved
  \begin{equation}
    \tr \omega_2^p  = p! \, x_0 \, |x|^{p-1} \, \tr \left( E_1 \cdot ... \cdot
    E_p \right) \, dx_1 \wedge ... \wedge dx_p 
  \label{G5-5.71}
  \end{equation}
  for all $x \in \R^{p+1} \setminus \{ 0 \}$.
 
 Next we calculate the second term in (\ref{G5-5.67}), where we again
  suppose $x' = (x_1,0,...,0)$:
  \begin{equation}
    \begin{split}
      &\tr \omega_1 \omega_2^{p-1}  |_{(x_0,x')} =
      \sum_{1\leq i_1 < ... < i_{p-1} \leq p} \: 
      \sum_{\sigma \in {\rm S}_{p-1}} \, \sgn \sigma
  \cdot \tr \big(  (x_0 - x_1 E_1 )^2 E_{i_{\sigma (1)}} 
      \cdot\\
  &\quad \cdot (x_0 - x_1 E_1 ) E_{i_{\sigma (2)}} \cdot 
      ... \cdot (x_0 - x_1 E_1 ) E_{i_{\sigma (p-1)}}  \big)
    \, dx_0 \wedge dx_{i_1} \wedge ... \wedge dx_{i_p}.
    \end{split}
  \end{equation}
  Now consider the trace terms on the right hand side of the last equation
  with $i_1=1 $ and $j = \sigma^{-1} (1)$:
  \begin{equation}
    \begin{split}
      \tr & \left(  (x_0 - x_1 E_1 )^2 E_{i_{\sigma (1)}} 
      \cdot (x_0 - x_1 E_1 ) E_{i_{\sigma (2)}} \cdot 
      ... \cdot (x_0 - x_1 E_1 ) E_{i_{\sigma (p-1)}}  \right)=\\
      & = \tr \left(  (x_0 + (-1)^j  x_1 E_1 ) \, (x_0 - x_1 E_1 )
      E_{i_{\sigma (j)}} \cdot ... \cdot (x_0 - x_1 E_1 ) 
      E_{i_{\sigma (p-1)}} \cdot \right. \\
      & \hspace{10mm} 
      \left. \cdot (x_0 - x_1 E_1 ) E_{i_{\sigma (1)}} \cdot ... \cdot 
      (x_0 - x_1 E_1 ) E_{i_{\sigma (j-1)}}  \right) \\
      & = |x|^{p-1} \tr \left( (x_0 + (-1)^j x_1 E_1 ) E_1 \cdot 
      E_{i_{\sigma(j+1)}} \cdot ...\cdot
      E_{i_{\sigma(p-1)}} E_{i_{\sigma(1)}}\cdot ... 
      \cdot E_{i_{\sigma(j-1)}} \right) = 0.
    \end{split}
  \end{equation}
  In case $i_1 > 1 $ an analogous argument proves
  \begin{equation}
    \begin{split}
      \tr & \left(  (x_0 - x_1 E_1 )^2 E_{i_{\sigma (1)}} 
      \cdot (x_0 - x_1 E_1 ) E_{i_{\sigma (2)}} \cdot 
      ... \cdot (x_0 - x_1 E_1 ) E_{i_{\sigma (p-1)}}  \right)=\\
      & = |x|^{p-1} \, \tr \left( (x_0 - x_1 E_1) \cdot E_{i_{\sigma (1)}}
      \cdot ... \cdot E_{i_{\sigma (p-1)}} \right) \\
      & = -x_1 \, |x|^{p-1} \, \tr \left( E_1 \cdot E_{i_{\sigma (1)}} 
      \cdot ... \cdot E_{i_{\sigma (p-1)}}\right) .
    \end{split}
  \end{equation}
  Hence by rotation symmetry
  \begin{equation}\begin{split}
    \tr \omega_1 \omega_2^{p-1}  = \,&|x|^{p-1} \, (p-1)! \, 
    \tr (E_1 \cdot ... \cdot E_p ) \, \\
   & dx_0 \wedge \sum_{j=1}^p
    \, (-1)^j \, x_j dx_1 \wedge ... \wedge \widehat{dx_j} \wedge ... dx_p  
                  \end{split}
  \label{G5-5.75}  
  \end{equation}
  holds. 
The assertion follows from  
\myref{G5-5.67}, \myref{G5-5.71} and \myref{G5-5.75}.
\end{proof}
From now on we choose the standard representation of 
$\C \ell_{2k-1}$ in $\C^{2^{k-1}}$. For this representation
one knows that 
\begin{equation}
  \tr ( E_1 \cdot ... \cdot E_{2k-1} ) = 2^{k-1} \, i^{-k} .
   \label{G5-5.76}
\end{equation}
This is part of the Berezin Lemma (cf.~\cite[Prop.~3.21]{BerGetVer:HKDO}),
but can also easily seen 
as follows. In the standard representation the complex volume 
element $i^k \, e_1 \cdot ... \cdot e_{2k-1}$ acts as identity.
Since the standard representation is of rank $2^{k-1}$ we obtain
\myref{G5-5.76}.
Now, by Proposition \ref{MPS1-4.1}
\begin{equation}
  \omega = |x|^{-p-1} \, p! \, 2^{k-1} \, i^{-k} \, \sum_{j=0}^p
  (-1)^j \, x_j \, dx_0 \wedge ... \wedge\widehat{dx_j}\wedge ... \wedge dx_p,
\end{equation}
and 
\begin{equation}
  \begin{split}
    \int_{S^{2k-1}} \, \omega & = \int_{B^{2k}} \, d (|x|^{p+1}\omega )= \\
    & = i^{-k} \, 2^{k-1} \, (2k-1)! \, \int_{B^{2k}}
    \, (p+1) \, d {\rm vol} \\
  %  & = i^{-k} \, 2^{k-1} \, (2k-1)! \, (2k) \, \frac{\pi^k}{\Gamma (k+1)} \\
    & = i^{-k} \, 2^k \, (2k-1)! \,  \frac{\pi^k}{(k-1)!} \\
    &=\frac{-1}{c_k},
  \end{split}
  \label{ML5-5.23}
\end{equation}
where $c_k$ was defined in Definition \plref{S5-5.1}.
%
%So we end up with the formula
%
%\begin{equation}
%  \frac{1}{(2 \pi i )^k} \, \frac{(k-1)!}{(2k-1)!} \, \int_{S^{2k-1}} \omega
%  = (-1)^k ,
%\end{equation}
%
%i.e.~$ c_k \int_{S^{2k-1}} \, \omega = -1$ with $c_k$ appropriate.

We note another consequence of our calculations.
Choose $a \in \R$ and let $f(x) = a + c(x)$, $x \in \R$.
Then $\left( f^* df\right)^p = \omega_2^p$ and
in view of \myref{G5-5.71} and \myref{G5-5.76}
\begin{equation}
\label{G4-5.5}
  \begin{split}
    \tr \left( (f^{-1} df)^p \right) & = (a^2 + |x|^2 )^{-p} \,
    \tr \left( (f^* df)^p \right) \\
    & = p! \, a (a^2 + |x|^2 )^{-\frac{p+1}{2}}
    \tr (E_1 \cdot ... \cdot E_p) \, dx_1 \wedge ... \wedge dx_p \\
    & = p! \, a (a^2 + |x|^2 )^{-\frac{p+1}{2}} \, 2^{k-1} \, i^{-k} \,
    dx_1 \wedge ... \wedge dx_p.
  \end{split}
\end{equation}
We are now able to calculate the integral of the $p$-form 
$\tr \left( (f^{-1} df)^p \right) $:
\begin{equation}
  \begin{split}
    \int_{|x| \leq R} \tr \left( (f^{-1} df)^p \right) & = 2^{k-1} \, i^{-k} \,
    p! \, a \int_{|x| \leq R} (a^2 + |x|^2 )^{-k} \, dx \\
    & = 2^{k-1} \, i^{-k} \, p! \, a \,
    \frac{(2k-1) \pi^{k - \frac 12}}{\Gamma (k + \frac 12)} \int_0^R    
    (a^2 + r^2 )^{-k} \, r^{2k-2}\, dr  \\
    & \underset{R \rightarrow \infty}{\longrightarrow} 
    2^{k-1} \, i^{-k} \, p! \, a \, 
    \frac{(2k-1) \pi^{k - \frac 12}}{\Gamma (k + \frac 12)} \int_0^\infty
    (a^2 + r^2)^{-k} \, r^{2k-2} \, dr.
  \end{split}
\end{equation}
Using the formula
\begin{equation}
  \int_0^\infty (1+ u^2)^\alpha u^\beta\, du = \frac12 \frac{\Gamma (- \alpha -
  \frac{\beta}{2} - \frac12) \,
  \Gamma (\frac{\beta +1}{2})}{\Gamma (- \alpha) }
\end{equation}
we find 
\begin{equation}
\begin{split}
  \int_{\R^p} \tr \left( (f^{-1} df)^p \right) & =
  2^{k-1} \, i^{-k} \, p! \, (2k-1) \, ( {\rm sgn} \, a) \, \frac 12 \,
  \frac{\pi^{k-\frac 12}}{\Gamma (k + \frac 12 )} \,
  \frac{\Gamma (\frac 12) \, \Gamma (k -  \frac 12 )}{\Gamma (k) } \\
  & = 2^{k-1} \, i^{-k} \, p! \, ({\rm sgn} \, a) \, \frac{\pi^k}{(k-1)!} 
   = \frac{- {\rm sgn} \, a}{2 c_k}. 
\end{split}
\label{G5-5.83}
\end{equation}
Now we are ready to relate 
the $\eta$-invariant to the spectral $\eta$-invariant
of an elliptic operator. 

First we briefly recall the regularized integral for functions
on $\R_+$ (cf.~\cite[Sec.~2.1]{Les:OFTCSAM}, \cite{LesTol:DODEBVP}).
Let $f:(0,\infty)\to\C$ be a locally integrable function having
$\log$--polyhomogeneous asymptotic expansions as $x\to 0$ and
as $x\to\infty$. Then one puts
\begin{equation}
    \regint_0^\infty f(x) dx :=\LIM_{a\to 0} \int_a^1 f(x) dx+
     \LIM_{b\to \infty} \int_1^\infty f(x) dx.
    \label{ML5-5.30}
\end{equation}
For such a function its ``Mellin transform''
\begin{equation}
      (\widetilde \cM f)(s):=\regint_0^\infty x^{s-1} f(x) dx
\end{equation}
is well--defined for $s\in \C$ and there is a discrete subset 
$A\subset \C$ such that $(\widetilde \cM f)\restr (\C\setminus A)$ extends
to a meromorphic function $\cM f$ on $\C$. 
For each $s$ one has
\begin{equation}
    (\widetilde \cM f)(s)={\rm Res}_0 (\cM f)(s),
   \label{ML5-5.32}
\end{equation}
where ${\rm Res}_0$ denotes the constant term in the Laurent expansion.
In particular, if $s$ is a regular point of $\cM f$ then
$(\widetilde \cM f)(s)=(\cM f)(s)$. 
We note that for $\ga\in \R, k\in\Z_+$
\begin{equation}
     \regint_0^\infty x^\ga \log^kx dx=0.
    \label{ML5-5.33}
\end{equation}
(see \cite[Sec.~2.1]{Les:OFTCSAM} for proofs of these facts).
Of course, there is a simple relation between the integral
\myref{ML5-5.30} and the integral \myref{ML5-4.3}.
Namely, if $f$ is an even $\log$--polyhomogeneous function on $\R$
then $2\reginttext\limits_0^\infty f(x)dx=\reginttext\limits_\R f(x) dx$.

\medskip
\comment{
Next let $D \in {\rm CL}^m (M,E)$ be a self-adjoint elliptic operator. For
$x> 0$ the operator $D^2+x$ is also elliptic and hence the $\zeta$--function
\begin{equation}
  s \mapsto \tr \left( D (D^2 + x)^{-s} \right) 
\end{equation}
extends  meromorphically to $\C$ with simple poles in
$s = \frac{\dim M + m - k}{2m}$, $k \in \Z_+$
(see e.g. \cite[Sec ???]{GruSee:WPPOAPSBP})\marginpar{???}. 

We put for $\alpha \in \C$
\begin{equation}
  {\rm tr}_{\zeta} \left( D(D^2 + x)^{-\alpha} \right) := 
  {\rm Res}_0 \tr \left. \left( D(D^2+x)^{-s} \right) \right|_{s=\alpha}.
\end{equation}
Obviously, we have for $\alpha \notin \Z_+$
\begin{equation}
  \begin{split}
    &\left( \frac{d}{dx} \right)^k {\rm tr}_{\zeta} \left( 
    D (D^2 + x)^{-\alpha} \right)
    = (- \alpha )\, (- \alpha -1) \cdot ... \cdot (-\alpha -k +1 ) \,
    {\rm tr}_{\zeta} \left( D (D^2 + x)^{-\alpha-k} \right) \\
    =& (- \alpha )\, (- \alpha -1) \cdot ... \cdot (-\alpha -k +1 ) 
    \, {\rm tr}_{{\rm L}^2} \left( D (D^2 + x)^{-\alpha-k} \right)
  \end{split}
\end{equation}
for $k$ large enough. 
Applying the resolvent expansion of an elliptic pseudodifferential
operator \cite[Sec. ???]{GruSee:WPPOAPSBP}\marginpar{???} we obtain
for $k\in\{1,2,3,...\}$ an asymptotic expansion 
\begin{equation}
  {\rm tr}_{\zeta} \left( D (D^2 +x)^{-k} \right) \sim_{x \rightarrow \infty}
  \sum_{{\rm Re}\, \beta \rightarrow - \infty} \,\sum_{j=0}^{k(\beta)}
  a_{\beta,j} \, x^\beta \, \log^j x,  
\end{equation}
where $k(\beta)\le 1$ for all $\beta$.
Thus the regularized integral introduced in \myref{ML5-5.30} can
be applied to the function $x\mapsto {\rm \tr}_\zeta
\left( D (D^2 +x)^{-k} \right)$ for $k\in \Z, k>1$.

\begin{lemma}
\label{S3-5.3}
  Let
  \begin{equation}
  \label{G4-5.2}
    \eta_D(s) := \sum_{\mu \in {\rm spec} \, D \setminus \{ 0 \} } \,
    {\rm sgn} \, \mu \: |\mu|^{-s}=\tr((D^2)^{-(s+1)/2})
  \end{equation}
  be the spectral $\eta$-function of $D$. Then we have for 
  $k \in \{ 1,2,3,...\}$ 
  \begin{equation}
  \label{G4-5.3}
    \eta (s) = \left( \begin{array}{c}k-1-\frac {s+1} 2 \\ k-1\end{array}
              \right)^{-1}
    \, \frac{ \sin \pi \frac {s+1} 2 }{\pi} \, \regint_0^{\infty} 
    \, x^{k-1-\frac {s+1} 2 } \, {\rm tr}_\xi \left( D (D^2 +x)^{-k} \right)
    \, dx,
  \end{equation}
  in particular
  \begin{equation}
  \label{G4-5.4}
    \eta (0) = \frac{\Gamma (k) }{\Gamma (k- \frac 12 ) \, \sqrt{\pi}}
    \regint_0^{\infty} x^{k-\frac 3 2} {\rm tr}_\xi \left( D ( D^2 +x)^{-k} \right)
    \, dx.
  \end{equation}
  Here $\regint$ denotes the regularized integral 
  (cf.~{\rm \cite[Sec.~2.1]{Les:OFTCSAM}, \cite{LesTol:DODEBVP}}).
\end{lemma}
\begin{proof}
  From
  \begin{equation}
    \begin{split}
      \lambda^{-z} & = \frac{ \sin\pi z}{\pi} \, \int_0^{\infty} \,
      x^{-z} \, (\lambda +x)^{-1} \, dx,  \quad 0 < {\rm Re} \, z < 1, \\
      \lambda^{-z} & = 
     %\left( \begin{array}{c} k-1-z \\ k-1\end{array}
     %    \right)^{-1}
     {\binom{k-1-z}{k-1}}^{-1}     
 \, \frac{\sin \pi z}{\pi} \, \int_0^{\infty} \,
      x^{k-1-z} \, (\lambda + x)^{-k} \, dx ,
      \quad 0 < {\rm Re} \, z < k,
    \end{split}
  \end{equation}
  and \myref{G4-5.2} we infer \myref{G4-5.3} for $k$ large enough
  and $s_0(D)<{\rm Re}\, s <2k-1$.
  Integration by parts gives \myref{G4-5.3} for all $k \in \{ 1,2,3,...\}$.
  Since 
$%\left( \begin{array}{c}k-1 \\ k-1- \frac {s+1} 2 \end{array}
 %\right)^{-1} 
   {\binom{k-1}{k-1-\frac{s+1}{2}}}^{-1}\, \frac{\sin \,
  \pi \frac {s+1} 2}{\pi} $ is regular and $\neq 0$ at $s =0$, and since
  $\eta (s)$ is regular at $s=0$ (cf.~\cite[Sec.~3.8]{Gil:ITHEASIT})
  we conclude that the meromorphic function
  \begin{equation}
    s \mapsto \regint_0^{\infty} \, x^{k-1- \frac {s+1} 2 }
    \, {\rm tr}_{\zeta} \left( D (D^2 + x)^{-k}  \right)  dx
  \end{equation}
  is regular at $s = 0$. This implies \myref{G4-5.4}.
\end{proof}
}

Next let $D:\cinf{E}\to \cinf{E}$ be a first order invertible
self-adjoint elliptic \emph{differential} operator. Then for
$k\in\Z_+, k>0,$ the family
\begin{equation}
  \Phi_{k,p}:\R^p\ni x \mapsto D (D^2 + |x|^2)^{-k} 
\end{equation}
lies in $\CL^{1-2k}(M,E,\R^p)$. Note that being a \emph{differential}
operator is really essential for this to be true.
We apply Theorem \plref{S1-3.4} to $\Phi_{k,p}$. Thus for each $k\in\Z_+, k>0,$
we have
\begin{equation}
       \TR(D(D^2+|\id_{\R^p}|^2)^{-k})=\TR(\Phi_{k,p})\in 
       \tilde\sym^{2k+\dim M}(\R^p)/\cP.
\end{equation}
We note for $p=1$ and $l>\frac{\dim M+1-2k}{2}$ the identity
\begin{equation}
   \left(\frac 1x \frac{\pl}{\pl x}\right)^l
    \TR(\Phi_{k,p})(x)=(-2)^l\frac{(k+l-1)!}{(k-1)!}
      \trltwo(D(D^2+|x|^2)^{-k-l}).
    \label{ML5-5.36}
\end{equation}
Note that in view of \myref{ML5-5.33} the Mellin transforms
$\widetilde \cM(\Phi_{k,1}), \cM(\Phi_{k,1})$, are well--defined.
For reasons of clarity we now write the various traces
on $\CL(M,E,\R^p)$ with a subscript, i.e.
$\TR_p, \overline{\rm Tr}_p, \widetilde{\rm Tr}_p$.
\begin{prop}
\label{S3-5.3}
  Let $D$ be as before and let
  \begin{equation}
  \label{G4-5.2}
    \eta_D(s) := \sum_{\mu \in {\rm spec} \, D } \,
    {\rm sgn} \, \mu \: |\mu|^{-s}=\tr((D^2)^{-(s+1)/2})
 \end{equation}
  be the spectral $\eta$-function of $D$. Then we have for 
  $k \in\Z_+, k>0,$ and ${\rm Re}\, s$ large
  \begin{equation}
  \label{G4-5.3}\begin{split}
    \eta (s) =& \left( \begin{array}{c}k-1-\frac {s+1} 2 \\ k-1\end{array}
              \right)^{-1}
    \, \frac{ \sin \pi \frac {s+1} 2 }{\pi} \, 2\regint_0^{\infty} 
    \, x^{2k-2-s} \, \TR_1\left( D (D^2 +\id_\R^2)^{-k} \right)(x) 
    \, dx,\\
     =& \left( \begin{array}{c}k-1-\frac {s+1} 2 \\ k-1\end{array}
              \right)^{-1}
    \, \frac{ \sin \pi \frac {s+1} 2 }{\pi} \, 2(\cM \Phi_{k,1})(2k-1-s).
                \end{split}
    \end{equation}
Furthermore,
  \begin{equation}
  \label{G4-5.4}\begin{split}
    \eta (0) = &\frac{2 \Gamma (k) }{\Gamma (k- \frac 12 ) \, \sqrt{\pi}}
    \regint_0^{\infty} x^{2k-2} \TR_1 \left( D ( D^2 +\id_\R^2)^{-k} \right)(x)
    \, dx\\
     =&\frac{ \Gamma (k) }{\Gamma (k- \frac 12 ) \, \sqrt{\pi}}
    \overline{\rm Tr}_1 (\id_\R^{2k-2}D(D^2+\id_\R^2)^{-k})\\
     =& \frac{\Gamma(k)}{\pi^k} \overline{\rm Tr}_{2k-1}(
        D(D^2+|\id_{\R^{2k-1}}|^2)^{-k}).
                \end{split}
  \end{equation}
\end{prop}
\begin{proof}
  From the identities
  \begin{equation}
    \begin{split}
      \lambda^{-z} & = \frac{ \sin\pi z}{\pi} \, \int_0^{\infty} \,
      x^{-z} \, (\lambda +x)^{-1} \, dx,  \quad 0 < {\rm Re} \, z < 1, \\
      \lambda^{-z} & = 
     %\left( \begin{array}{c} k-1-z \\ k-1\end{array}
     %    \right)^{-1}
     {\binom{k-1-z}{k-1}}^{-1}     
 \, \frac{\sin \pi z}{\pi} \, \int_0^{\infty} \,
      x^{k-1-z} \, (\lambda + x)^{-k} \, dx ,
      \quad 0 < {\rm Re} \, z < k,
    \end{split}
  \end{equation}
  and \myref{G4-5.2} we infer \myref{G4-5.3} for $k$ large enough
  and $s_0(D)<{\rm Re}\, s <2k-1$.
  Integration by parts gives \myref{G4-5.3} for all $k \in \{ 1,2,3,...\}$
 (cf. \myref{ML5-5.36}).
  Since 
$%\left( \begin{array}{c}k-1 \\ k-1- \frac {s+1} 2 \end{array}
 %\right)^{-1} 
   {\binom{k-1-\frac{s+1}{2}}{k-1}}^{-1}\, \frac{\sin \,
  \pi \frac {s+1} 2}{\pi} $ is regular and $\neq 0$ at $s =0$, and since
  $\eta (s)$ is regular at $s=0$ (cf.~\cite[Sec.~3.8]{Gil:ITHEASIT})
  we conclude from \myref{G4-5.3} that the meromorphic function
  $(\cM\Phi_{k,1})(2k-1-s)$ 
  is regular at $s = 0$ and thus in view of \myref{ML5-5.32}
  we arrive at the first equality of \myref{G4-5.4}.
  The second equality of \myref{G4-5.4} is trivial. To prove the
third one we note that from the uniqueness statement and
(iii) of Theorem \plref{S1-3.4} we conclude the identity
\begin{equation}
    \Phi_{k,p}(x)=\Phi_{k,1}(|x|)
\end{equation}
and thus
\begin{equation}
   \begin{split}
     &\frac{\Gamma(k)}{\pi^k} \overline{\rm Tr}_{2k-1}(
        D(D^2+|\id_{\R^{2k-1}}|^2)^{-k})\\
     =& \frac{\Gamma(k)}{\pi^k} \regint_{\R^{2k-1}}\TR_{2k-1}
         (D(D^2+|\id_{\R^{2k-1}}|^2)^{-k})(x) dx\\
     =&\frac{\Gamma(k)}{\pi^k} \LIM_{R\to\infty}\int_{|x|\le R}\TR_{2k-1}
         (D(D^2+|\id_{\R^{2k-1}}|^2)^{-k})(x) dx\\
    = & \frac{\Gamma(k)}{\pi^k}\LIM_{R\to\infty}
          \frac{(2k-1) \pi^{k-1/2}}{\Gamma(k+\frac 12)}
          \int_0^R \TR_1(D(D^2+\id_\R^2)^{-k})(x)x^{2k-2}dx\\
    =&\frac{ \Gamma (k) }{\Gamma (k- \frac 12 ) \, \sqrt{\pi}}
    \overline{\rm Tr}_1 (\id_\R^{2k-2}D(D^2+\id_\R^2)^{-k}).
   \end{split}
\end{equation}
\end{proof}
\begin{proposition}\label{ML5-5.6}
  Let $D:\cinf{E}\to\cinf{E}$ be a first order invertible 
  self-adjoint elliptic 
  differential operator. Let $c: \R^{2k-1} \rightarrow {\rm M} (2^{k-1}, \C)$
  be the standard Clifford representation (see \myref{G5-5.76}).  
  Then the family
  $\cD_\pm(\mu) := D \pm c (\mu)$ lies in ${\rm CL}^1 (M,E, \R^{2k-1} )$ 
  and we have
  \begin{equation}
    \eta (\cD_\pm) = \mp \eta (D).
  \end{equation}
\end{proposition}
\begin{remark}
\mbox{} 
\begin{enumerate}\renewcommand{\labelenumi}{{\rm \arabic{enumi}.}}
\item
  This generalizes \cite[Prop.~5]{Mel:EIFPO}.
  Even for $k =1$ it is more general than there.
  In contrast to \cite{Mel:EIFPO} 
  our proof does not use the local index theorem,
  hence is not restricted to Dirac operators. 
\item Since the spectral $\eta$--invariant is a regularization
of the non--convergent sum
$$\sum_{\mu \in {\rm spec} \, D \setminus \{ 0 \}} {\rm sgn}\, \mu,$$
this result is formally a consequence of the integral formula
\myref{G5-5.83}.
\end{enumerate}
\end{remark}
\begin{proof}
  Since $D$ is a differential
  operator we have $\cD \in {\rm CL}^1 (M,E, \R^{2k-1} )$, 
  where  being differential is really essential here.
  Note that the complete symbol $\sigma_\cD(x,\xi,\mu)$ is
  (affine) linear in $(\xi,\mu)$.
  From \myref{G4-5.5} and the previous proposition we conclude
  \begin{equation}
    \begin{split}
      \eta (\cD_\pm) & = 2 \, c_k \, \regint_{\R^{2k-1}} \, \TR_{2k-1}
        \left( \left( (\pm D +
          c(\cdot) )^{-1} \, d c(\cdot) \right)^{2k-1} \right) (\mu) \, d \mu \\
      & = \pm 2^k \, i^{-k} \, c_k \, (2k-1) ! \,
      \regint_{\R^{2k-1}} \TR_{2k-1}(D(D^2+|\id_{\R^{2k-1}}|^2)^{-k})(x) dx\\
      & = \ \pm 2^k \, i^{-k} \, c_k \, (2k-1) ! \,
      \frac{\pi^k}{\Gamma(k)} \, \eta (D) \\
      & = \mp \eta (D),
    \end{split}
  \end{equation}
  and we are done.
\end{proof}

\begin{remark}
(i)
$\eta_1$ is an additive homomorphism from the group of invertible elements of 
${\rm CL}^*(M,E,\R)$ into $\C$ \cite[Prop. 4]{Mel:EIFPO}.
This follows immediately from the conjugation invariance of
the trace $\overline{\Tr}$. Namely, given invertible
$A,B\in {\rm CL}^*(M,E,\R)$ then
\begin{equation}
\begin{split}
    (AB)^{-1} d(AB)&= B^{-1}(A^{-1} dA) B+ B^{-1} dB\\
                   &=: \go_1+\go_2.
\end{split} \label{ML5-5.45}
\end{equation}
Thus
\begin{equation}\begin{split}
     \eta_1(AB)&=  \frac{1}{\pi i} 
               \overline{\Tr}( B^{-1}(A^{-1} dA) B+ B^{-1} dB)\\
    &= \eta_1(A)+\eta_1(B).
                \end{split}
\end{equation}

However, $\eta_k$ is not additive for $k\ge 2$.
We illustrate this in the case $k=2$:
the $1$--forms $\go_1, \go_2$ of \myref{ML5-5.45} have the following
properties:
\begin{equation}
\begin{split}
    &d\go_1=-\go_1^2-\go_1\wedge\go_2-\go_2\wedge\go_1,\quad
    d\go_2=-\go_2^2,\\
    &d(\go_1+\go_2)=-(\go_1+\go_2)^2,\quad
    d(\go_1\wedge \go_2)=-(\go_1+\go_2)\wedge\go_1\wedge\go_2.
\end{split}
\end{equation}
Consequently
\begin{equation}
\begin{split}
     \TR((\go_1+\go_2)^3)-\TR(\go_1^3)-\TR(\go_2^3)&=
       3\TR((\go_1+\go_2)\wedge\go_1\wedge\go_2)\\
       &= -3 d\TR(\go_1\wedge\go_2)
\end{split}
\end{equation}
and 
\begin{equation}
  \eta_2(AB)  = \eta_2(A)  + \eta_2 (B) -6c_2 
      \widetilde{\TR} (\omega_1 \wedge\omega_2). 
\end{equation}  
  So the defect of the additivity of $\eta$ is a symbolic term.

\bigskip (ii) Finally we add some remarks concerning the divisor flow
(cf.~\cite[Sec. 9]{Mel:EIFPO}). Following \cite[Sec.~8]{Mel:EIFPO}, 
the right hand side of the variation formula Proposition \plref{ML5-5.3}
can be defined if $A_s$ is elliptic but not necessarily invertible.
Namely, choose a smooth family of parametrices 
$Q_s\in\CL^{-m}(M,E,\R^{2k-1})$, i.e.
$$Q_sA_s-I,\quad A_sQ_s-I\in \CL^{-\infty}(M,E,\R^{2k-1}),$$
and put
\begin{equation}
     {\rm v\eta}(A_s):= 2(2k-1) c_k \widetilde{\Tr} \left(
    (Q_s \partial_s A_s) (Q_s d A_s)^{2k-2} \right).       
 \end{equation}
In view of Proposition \plref{S3-4.6}   ${\rm v\eta}(A_s)$
is independent of the choice of $Q_s$.

Now fix an elliptic and invertible $A\in\CL(M,E,\R^{2k-1})$. If
$B\in{\rm Ell}_m(A)$, the component of $A$ in the elliptic elements
of order $m$, then according to \cite[Def.~5]{Mel:EIFPO}
the divisor flow was defined to be
\begin{equation}
   {\rm DF}(B,A):=\frac 12\left( \eta_k(B)-\eta_k(A)-\int_0^1{\rm v\eta}(B_s)ds
   \right)
\end{equation}
where $B_s$ is any smooth family in ${\rm Ell}_m(A)$
with $B_0=A, B_1=B$.

However, it is not immediately clear, whether ${\rm DF}(B,A)$ is independent
of the particular choice of a path $B_s$. 
There seems to be some evidence that this might not be the case.
Let us give an example for $M$ a point. Although this is quite an exceptional
case it at least shows where the problems are:

Let $f\in\cinf{\R}$ with $f'\in\cinfz{\R}$, and
$\lim\limits_{\gl \to\pm\infty} f(\gl) \neq 0$. If $f$ is invertible, then
\begin{equation}
       \eta(f):=\frac{1}{\pi i} \int_\R \frac{f'}{f}(\gl) d\gl.
\end{equation}
If $f_s$ is a family of such functions,  the variation
formula reads
\begin{equation}
      \begin{split}
    \pl_s\eta(f_s)=& \frac{1}{\pi i} \int_\R \pl_\gl \left( 
    f_s^{-1} \pl_s f_s \right) \, d\gl \\
     =& \frac{1}{\pi i} \left(\frac{\pl_s f_s}{f_s}(+\infty)-
          \frac{\pl_s f_s}{f_s}(-\infty)\right)\\
     =:& {\rm v\eta}(f_s).
      \end{split}
\end{equation}
Now let 
\begin{equation}
     f_s(\gl):=\left\{\begin{array}{ll}
           1,&\gl\ge 1,\\
        e^{2\pi i s},& \gl\le s,\\
      e^{2 \pi i \gl}, & s\le \gl\le 1.
                    \end{array}\right.
\end{equation}
The two points at which $f_s$ is just continuous but not smooth can
easily be smoothed out.

Then one calculates
\begin{equation}
     \int_0^1 {\rm v\eta}(f_s) \, ds =-2.
\end{equation}
On the other hand 
\begin{equation}
    g_s := (1-s)f_0+s f_1,\quad 0\le s\le 1,
\end{equation}
also is an elliptic family joining $f_0$ and $f_1$. Here, elliptic
means invertible outside a compact set. However,
\begin{equation}
       \int_0^1 {\rm v\eta}(g_s)\, ds=
       0\neq \int_0^1{\rm v\eta}(f_s)\, ds
\end{equation}
proving the path dependence of the divisor flow.
\end{remark} 

%%% Local Variables: 
%%% mode: latex
%%% TeX-master: "LesPflTAPDPOEI"
%%% End: 

%
\addcontentsline{toc}{section}{References}
\newcommand{\und}{{\rm and }}

%%% Local Variables: 
%%% mode: latex
%%% TeX-master: "LesPflTAPDPOEI"
%%% End: 

\bigskip\noindent
Humboldt Universit\"at zu Berlin\\
Institut f\"ur Mathematik\\ 
Unter den Linden 6\\
D--10099 Berlin
\bigskip\noindent
\begin{tabbing}
Email--addresses: \=
   lesch@mathematik.hu-berlin.de\\
\>pflaum@mathematik.hu-berlin.de\\[1em]
Internet: \>  http://spectrum.mathematik.hu-berlin.de/$\sim$lesch\\
 \>http://spectrum.mathematik.hu-berlin.de/$\sim$pflaum
\end{tabbing}
\end{document}